\documentclass[11pt,leqno]{amsart}
\usepackage[T1]{fontenc}
\usepackage{ae}
\usepackage{amsmath}
\usepackage{amssymb}
\usepackage{chicago}
\usepackage{latexsym,graphicx}
\usepackage{graphicx,multirow,rotating}
\renewcommand{\baselinestretch}{2.0}

\renewcommand{\baselinestretch}{1}
\newtheorem{theo}{Theorem}[section]
\newtheorem{lem}{ Lemma}[section]

\newtheorem{exem}{Example}
\newcommand{\eref}[1]{(\ref{#1})}
\newcommand\ind{{{{1}}\hspace{-1,1mm}{\mathrm I}}}
\newcounter{hypc}
\newcommand{\hypothese}[1]{\stepcounter{hypc}\tag{$\mathbf{A_{\thehypc}}$}
\label{#1}}

\newcounter{condc}
\newcommand{\condition}[1]{\stepcounter{condc}\tag{$\mathbf{C_{\thecondc}}$}
\label{#1}}

\def\build#1_#2^#3{\mathrel{\mathop{\kern 0pt#1}\limits_{#2}^{#3}}}
\def\cvl{\build{\ \longrightarrow\ }_{\nti}^{{\mathcal L}}}
\def\cvps{\build{\ \longrightarrow\ }_{\nti}^{{P.S.}}}
\def\cvp{\build{\ \longrightarrow\ }_{\nti}^{{\mathbb{P}}}}
\def\cvL1{\build{\ \longrightarrow\ }_{\nti}^{{\mathbb{L}^1}}}
\def\cv{\build{\ \longrightarrow\ }_{\nti}^{}}
\def\nti{n \rightarrow \infty}

\begin{document}
\title[Hazard function estimation with mismeasured covariate]{Semi-parametric estimation of the hazard function in a
  model with covariate measurement error}
\author{Marie-Laure Martin-Magniette$^{1,2}$, Marie-Luce TAUPIN$^{3}$}
\maketitle
\footnote*{\noindent I.N.A. Paris-Grignon, Math\'ematique et Informatique
Appliqu\'ees, Paris, \newline  $^2$ I.N.R.A.%Institut National de la Recherche
%Agronomique
, Unit\'e
de Recherche en Genomique V\'eg\'etale, Evry, France, \newline
$^3$Universit\'e Paris-Sud, Orsay and Universit\'e Paris Ren\'e
  Descartes, IUT de Paris}
%\begin{center}
%\today
%\end{center}

\begin{abstract}
We consider a model where  the failure hazard function, 
conditional on a time-independent covariate $Z$
is given by $R(t,\theta^0|Z)=\eta_{\gamma^0}(t)f_{\beta^0}(Z)$, with
$\theta^0=(\beta^0,\gamma^0)^\top\in \mathbb{R}^{m+p}$. The baseline hazard
function $\eta_{\gamma^0}$ and relative risk $f_{\beta^0}$ belong both to  parametric
families. The covariate $Z$ is measured with an error through an
additive error model $U=Z+\varepsilon$ where $\varepsilon$ is a random variable, independent from $Z$, with
known density $f_\varepsilon$.
We observe a $n$-sample $(X_i, D_i, U_i)$, $i=1,\ldots,n$, where $X_i$ is 
the minimum between the failure time and the censoring time, 
and $D_i$ is the censoring indicator.
We aim at estimating $\theta^0$
in presence of the unknown density $g$ of the covariate $Z$ using the observations
$(X_i, D_i, U_i)$, $i=1,\ldots,n$. Our estimation procedure 
based on least squares criterion provide two estimators of $\theta^0$.
The first one is obtained by minimizing an estimation of the least squares criterion where
$g$ is estimated by density deconvolution. We give upper bounds for  its risk that
depend on the smoothness properties of $f_\varepsilon$ and $f_\beta(z)$ as a function of $z$. We derive from this construction sufficient conditions
that ensure the $\sqrt{n}$-consistency. 
The second estimator is constructed under conditions ensuring that the least squares
criterion can be directly estimated with the parametric rate.
We propose a deep study of examples considering 
various type of relative risks $f_{\beta}$ and various types of error
density $f_\varepsilon$. We show in particular that in the Cox model and in
the excess
risk model, the estimators are
$\sqrt{n}$-consistent and asymptotically Gaussian estimators of $\theta^0$
whatever is $f_\varepsilon$.
\end{abstract}

\begin{abstract}
Consid\'erons un mod\`ele \`a risque instantan\'e 
mod\'elis\'e par la relation
$R(t,\theta^0|Z)=\eta_{\gamma^0}(t)f_{\beta^0}(Z)$, o\`u
$\theta^0=(\beta^0,\gamma^0)^\top\in
\mathbb{R}^{m+p}$. Le risque de base 
$\eta_{\gamma^0}$ et la fonction de risque relatif $f_{\beta^0}$ appartiennent  \`a des 
 familles param\'etriques.
La covariable $Z$ est mesur\'ee avec une erreur au travers de la relation $U=Z+\varepsilon$, $\varepsilon$ 
\'etant une variable al\'eatoire, ind\'ependante de $Z$, de densit\'e connue $f_\varepsilon$. 
Nous disposons d'un $n$-\'echantillon $(X_i,D_i,U_i)$, $i=1,\ldots, n$ o\`u $X_i$ est le
minimum entre le temps de survie et le temps de censure
  et $D_i$ est l'indicateur de censure. 
Notre but est d'estimer $\theta^0$, en pr\'esence la
densit\'e inconnue $g$, de la covariable $Z$, en utilisant les observations
$(X_i,D_i,U_i)$, $i=1,\ldots, n$.
Notre m\'ethode d'estimation, fond\'ee sur le crit\`ere des moindres carr\'es
nous fournit deux estimateurs. Pour le premier, nous \'etablissons des bornes
sup\'erieures du risque d\'ependant des r\'egularit\'es de la densit\'e des
erreurs $f_\varepsilon$ et de la fonction de risque relatif, comme
fonction de $z$.  Nous en d\'eduisons des conditions suffisantes pour
atteindre la vitesse param\'etrique. Le deuxi\`eme estimateur est construit
sous des hypoth\`eses assurant que le crit\`ere des moindres carr\'es
peut \^etre estim\'e \`a la vitesse param\'etrique.
Au travers d'exemples, nous
\'etudions les propri\'et\'es des estimateurs ainsi que les
conditions assurant la $\sqrt{n}$-consistance 
pour des fonctions de risque relatif
et des densit\'es d'erreurs vari\'ees. 
En particulier, dans le mod\`ele de Cox et dans le mod\`ele 
d'exc\`es de
risque, les estimateurs construits
sont  $\sqrt{n}$-consistants et asymptotiquement gaussiens, quelle que 
soit la loi des erreurs $\varepsilon$.
\end{abstract}

nonparametric estimation, excess risk model, Cox model
censoring, survival analysis.

\noindent {\it MSC Classifications (2000): Primary 62G05, 62F12,62N01, 62N02; 
Secondary   62J02.}

\bigskip

\bibliographystyle{chicago}

\date{}
\def\Cov{\mathop{\rm Cov}\nolimits}% % A sample of proper declaration of
\def\Argmax{\mathop{\rm Arg\,max}\limits}%    a math operator.
\def\Argmin{\mathop{\rm Arg\,min}\limits}%
\def\Var{\mathop{\rm Var}\nolimits}%
\newcommand {\LL}{\mbox{I\negthinspace L}}
\renewcommand{\baselinestretch}{1.1}
%\renewcommand{\baselinestretch}{2.0}
%\input{/home/mlm/DOC/MLUCE/com}
%\input{com}
%\allowdisplaybreaks
%\maketitle

\newpage

\section{Introduction}

In a proportional hazard model the hazard function is defined by
\begin{eqnarray}\label{model}
R(t,\theta^0| Z)=\eta_{\gamma^0}(t)f_{\beta^0}(Z),
\end{eqnarray}
where $\eta_{\gamma^0}$ is the baseline hazard function and  
$f_{\beta^0}$ is the relative risk, \textit{i.e.} the risk associated with the
value of the covariate $Z$ and relative to the risk under standard
condition given by $f_{\beta^0}(0)=1$. In this paper we consider general
relative risk $f_{\beta}$ with a special interest in
$f_\beta(z)=\exp(\beta z)$ and $f_\beta(z)=1+\beta z$ which define
respectively the Cox model and the model of excess relative risk. 
The functions $\eta_{\gamma^0}$ and $f_{\beta^0}$ belong both to 
parametric families and $\theta^0=(\beta^0,\gamma^0)^\top$ belongs to
the interior of a compact set $\Theta=\mathbb{B}\times \Gamma\subset
\mathbb{R}^{m+p}$. To ensure that the hazard function is a positive function, 
we assume that both are positive functions.

We are interested in the estimation of  $\theta^0$ when the covariate $Z$ is
measured with error. If $Z$ were measured without error, we
would consider a cohort of $n$ individuals during a fixed time interval $[0, \tau]$.
 For each individual, 
we would observe a triplet $(X_i, D_i, Z_i)$, where 
$X_i=\min(T_i,C_i)$ is the minimum between the failure time $T_i$ 
and the censoring time $C_i$, $D_i=\ind_{T_i \leq C_i}$ denotes 
the failure indicator, and $Z_i$ is the value of the covariate.
In this paper we consider that the covariate is mismeasured. 
For example the covariate $Z$ is a stage of a disease, not correctly diagnosed or 
a dose of ingested pathogenic agent, not correctly evaluated, so that the 
error range between the unknown dose and the evaluated dose  is sizeable. 
In this context, the available observation for each individual is the triplet
$\Delta_i=(X_i,D_i,U_i)$ where $U_i$ is an evaluation of the unobservable
covariate $Z_i$. The random variables $U$ and $Z$ are related by the error
model defined by
\begin{eqnarray}
\label{erreur}
U=Z+\varepsilon,
\end{eqnarray}
where $\varepsilon$ is a centered random variable, independent of
$Z$, $T$, and $C$. The density of $\varepsilon$ is known and denoted by $f_\varepsilon$.
Our aim is thus to estimate the parameter $\theta^0=(\beta^0,\gamma^0)^\top$ from the $n$-sample
$(\Delta_1,\ldots,\Delta_n)$ 
 in the presence of the unknown density $g$ of the unobservable covariate
$Z$, seen as a nuisance parameter  belonging to a 
functional space.

%In this setting we assume that the failure time  $T$ and the censoring time $C$, 
%conditional on $ Z$  and  $U$, are independent,the distribution of the failure
%time $C$, conditional on $Z$ and $U$, does not depend on $Z$ and  $U$ the distribution of the failure time $T$ conditional on $ Z$
%and $U$ does not  depend on $U$.

\subsection{Previous known results and ideas}

Models with measurement errors are deeply studied since the $50$'s with the
first papers of Kiefer and Wolfowitz \citeyear{KIEFWOLF} and Reiers{\o}l
\citeyear{REIER} for regression models with errors-in-variables.  We refer to Fuller \citeyear{Fuller} and  Carroll \textit{et al.}
\citeyear{Carrolletal95} for a presentation of such models and results related
to measurement error models. The interest for
survival models when covariates are subject to measurement errors is
more recent. 

To take into account that the
covariate $Z$ is measured with error, the first idea is 
simply to replace $Z$ with the observation $U$ in the score function defined
 by
\begin{eqnarray}
\label{vraispart}
L_n^{(1)}(\beta,Z^{(n)})=\frac{1}{n}\sum_{i=1}^n\int_0^\tau\left(
  \frac{f_{\beta}^{(1)}(Z_i)}{f_{\beta}(Z_i)}-\frac{
    \sum_{j=1}^nY_j(t)f_{\beta}^{(1)}(Z_j)
}{\sum_{j=1}^n
 Y_j(t)f_{\beta}(Z_j)}\right)dN_i(t), 
\end{eqnarray}
where $N_i(t)=\ind_{X_i \leq t, D_i=1}$, $Y_i(t)=\ind_{X_i \geq t}$,
$Z^{(n)}=(Z_1,\ldots,Z_n)$, and where $f^{(1)}_\beta$ is the first derivative
of $f_\beta$ with respect to $\beta$. We refer to Gill and Andersen
\citeyear{ANDERSENGILL} for futher details on \eref{vraispart}.
This method, named the naive method, is known to provide, even in the Cox
model, a
biased estimator of $\beta^0$. This comes from 
the fact that $$\lim_{\nti}\mathbb{E}[L_n^{(1)}(\beta^0,U^{(n)})]\not
=\lim_{\nti}\mathbb{E}[L_n^{(1)}(\beta^0,Z^{(n)})]=0.$$
To our knowledge, all
previously known results about consistency for
the semi-parametric estimation of the hazard function when the covariate is mismeasured are 
obtained in the Cox model. Let us present those results.
Various authors propose estimation
procedures based on corrections of the score function $L_n^{(1)}(\beta, U^{(n)})$. 
Among them, one can cite 
%Hu \textit{ et al.} \citeyear{HuTsiatisDavidian98}
%who consider this problem under additionnal assumptions (prior knowledge,
%parametric assumptions,...) or 
Kong \citeyear{KONG} who calculates the asymptotic bias of the
naive estimator obtained by
minimization of $L_n^{(1)}(\beta^0,U^{(n)})$, and defines an adjusted estimator. His  estimator  is not 
consistent, but a simulation study indicates that it is less biased 
than the  naive estimator. In the same context, Buzas \citeyear{BUZAS}
proposes an unbiased score function, and shows throughout a simulation study
that his method yields to an estimator
with a small bias. 
Following the approach developped first by Stefanski \citeyear{STEFANSKI89} and 
Nakamura \citeyear{NAKAMURA90} for generalized linear models,  Nakamura
\citeyear{NAKAMURA92} constructs 
an approximately  corrected partial score likelihood, defined by
$L_n^{(1)}(\beta,U^{(n)})+\sigma^2 \beta N(\tau), $
where $N(\tau)$ is the number of failures in the interval $[0, \tau]$ and
where $\varepsilon$ is a centered Gaussian random variable with variance $\sigma^2.$
Under the error model defined in \eref{erreur}, this correction is
based  on the facts : 
\begin{eqnarray}
\lim_{\nti}\!\!\!\!\!\!&&\!\!\!\!\!\!\mathbb{E}[L_n^{(1)}(\beta,Z^{(n)})]\mbox{ only depends on }\mathbb{E}(Z)
\mbox{ and }\mathbb{E}[\exp(\beta Z)]\label{cox1},\\
\mathbb{E}(Z)&=&\mathbb{E}(U)\label{cox2}\\ 
\mathbb{E}[\exp(\beta U)]&=&\mathbb{E}[\exp(\beta Z)]\mathbb{E}[\exp(\beta
\varepsilon)]\label{cox3}.\end{eqnarray}
Kong and Gu \citeyear{KONGGU} prove that the Nakamura
\citeyear{NAKAMURA92}'s  estimator  is a $\sqrt{n}$-consistent and
asymptotically Gaussian estimator of $\beta^0$. One can also cite Augustin \citeyear{Augustin04} who proposes an exact
correction of the log-likelihood function.

Again in the Cox model, an extension of the previously mentionned works
is presented in Hu and
Lin \citeyear{HuLin02}. They obtain a broad class of consistent estimators for the regression 
parameter when $U$ is measured on all study individuals
and the true covariate is ascertained on a randomly selected validation set.
A nonparametric correction approach  of the partial score function is also developped by Huang and Wang
\citeyear{HuangWang00} when repetitions are available.

We point out that those results 
strongly depend on the exponential form of the relative risk of the Cox
model, through the
use of \eref{cox1}-\eref{cox3} and the
extension of such methods to other relative risks is not concluding. 
For instance, in the model of excess relative risk without errors, the hazard
function is defined by
$R(t,\theta^0|Z)=\eta_{\gamma^0}(t)(1+\beta^0 Z)$ and the score function is
given by
\begin{eqnarray*}
L_n^{(1)}(\beta,Z^{(n)})=\frac{1}{n}\sum_{i=1}^n\int_0^\tau 
\left(\frac{Z_i}{1+\beta Z_i}-
\frac{ \sum_{i=j}^nY_j(t)Z_j}{\sum_{j=1}^n Y_j(t)(1+\beta Z_j)}\right)dN_i(t).
\end{eqnarray*} 
In this model, the naive method also provides biased estimator of $\beta^0$, since
$$\lim_{\nti}\mathbb{E}[L_n^{(1)}(\beta^0, U^{(n)})] \neq
\lim_{\nti}\mathbb{E}[L_n^{(1)}(\beta^0, Z^{(n)})]=0.$$
Indeed, easy calculations combined with the Law of Large Numbers give that
the limit $\lim_{\nti}\mathbb{E}[L_n^{(1)}(\beta, Z^{(n)}) ]$ 
depends on  $\mathbb{E}[Z/(1+\beta Z)]$ whereas the limit 
$\lim_{\nti}\mathbb{E}[L_n^{(1)}(\beta,U^{(n)}) ]$ depends on  $\mathbb{E}[U/(1+\beta U)]$.
Since the error model \eref{erreur} does not provide any expression
of $\mathbb{E}[Z/(1+\beta Z)]$  related to $\mathbb{E}[U/(1+\beta U)]$,
corrections analogous to the ones proposed in the Cox model cannot be
exhibited. In other words, it seems impossible to find a function $\Psi(\beta,U)$, independent of the
unknown density $g$ satisfying that $\mathbb{E}(\Psi(\beta,U))=\mathbb{E}[Z/(1+\beta Z)]$.
Consequently the methods proposed in the Cox model, by
Nakamura \citeyear{NAKAMURA92}, Kong and Gu \citeyear{KONGGU}, Buzas
\citeyear{BUZAS}, Lin \citeyear{HuLin02}, Huang and Wang
\citeyear{HuangWang00} or by Augustin \citeyear{Augustin04}
cannot
be applied to the model of excess relative risk and \textit{a fortiori} to a
model with a general
relative risk.
As a conclusion, methods based on a
correction of the partial score likelihood
\eref{vraispart} where $Z$ is replaced with $U$ seem
not concluding in a general setting.

An other possible way is to consider a partial log-likelihood related 
to the filtration generated by 
the observations.
If the covariate 
$Z$ were observable, then the filtration at time $t$, generated by the observations 
would be $\sigma\{Z, N(s), \ind_{X > s}, 0 \leq s \leq t\leq \tau\},$ 
and the intensity of the censored process $N(t)$ with respect to this
filtration would equal
$\lambda(t,\theta^0,Z)=\eta_{\gamma^0}(t)Y(t)f_{\beta^0}(Z).$
In case of covariate measurement error,  $Z$ is unobservable and 
only  the evaluation $U$ is available. In this context, 
the filtration generated by the observations is 
$\mathcal{E}_t=\sigma\{U, N(s), \ind_{X>s}, 0 \leq s \leq t\leq \tau\},$
%included in the following filtration
%$$\mathcal{F}_t=\sigma\{Z, U, N(s), \ind_{X>s}, 0 \leq s \leq t \leq \tau\}.$$
and the intensity of the censored process $N(t)$ with respect to the filtration  $\mathcal{E}_t$
equals
\begin{eqnarray*}
\mathbb{E}[\lambda(t,\theta^0,Z)|\mathcal{E}_t]=\eta_{\gamma^0}(t)Y(t)\mathbb{E}[f_{\beta^0}(Z)|\sigma(U,\ind_{T\geq t})].
\end{eqnarray*}
This is studied in the Cox model by Prentice \citeyear{PRENTICE82} who proposes the following  induced failure hazard 
function 
\begin{equation*}\label{induite}
\eta_{\gamma^0}(t)\mathbb{E}[\exp(\beta^0 Z)|\sigma(\ind_{T\geq t}, U)].
\end{equation*}
The presence of  $\{T \geq t\}$ in the conditioning usually 
implies that the induced partial log-likelihood  has not explicit 
form. When the marginal distribution of $Z$ given
$\{T\geq t, U\}$ is specified at each time $t$, Prentice \citeyear{PRENTICE82} proposes an approximation of the induced partial 
log-likelihood independent of the baseline hazard function. Nevertheless
this approximation is appropriate only when the disease is rare. 
Tsiatis \textit{et al.} \citeyear{TSIATIS} propose  another approximation of the 
induced partial log-likelihood, but emphasise that their method  cannot 
be used for the model of excess relative risk. 

In the Cox model with missing covariate Pons \citeyear{PONS} uses also
the partial likelihood. She proposes an estimator based on
\eref{vraispart} where $\mathbb{E}[\exp(\beta Z)|\sigma(U,\ind_{T\geq t})]$ is replaced with
$\mathbb{E}[\exp(\beta Z)|U]$ by considering that
\begin{eqnarray}
\label{pons}
\mathbb{E}[\exp(\beta Z)|\sigma(U,\ind_{T\geq t})]=\mathbb{E}[\exp(\beta Z)|U].\end{eqnarray}
Nevertheless, obvious examples can be exhibited to prove that the equality
\eref{pons} does not hold in a general setting.
As a conclusion, the partial likelihood related to
the filtration $\mathcal{E}_t$ seems unusable since it
is difficult to separate the estimation of $\beta^0$ from the estimation of
$\eta_{\gamma^0}$.

\subsection{Our results}

%%%%%%%%%%%%%%%%%%%%%%%%%%%%%%%%%%%%%%%%%%%%%%%%%%%%%%%%%%%%%%%%%%%%%%%%%%%%%%%%%%%%%%%%%%

Our estimation procedure is based on the estimation of least squares criterion
using deconvolution methods. More precisely, using the
observations $\Delta_i=(X_i,D_i,U_i)$ for $i=1,\ldots,n$, we estimate the least squares criterion 
\begin{eqnarray}
S_{\theta^0,g}(\theta)&=&
\mathbb{E}\left(f^2_{\beta}(Z)W(Z) \int_0^\tau Y(t)\eta^2_{\gamma}(t)dt\right)-2
\mathbb{E}\left(f_{\beta}(Z)W(Z)\int_0^\tau \eta_{\gamma}(t)dN(t)\right)\label{S0}.
\end{eqnarray}
The function $W$ is a positive weight function to be suitably chosen such that $Wf_\beta$
and its derivatives up to order 3 with respect to $\beta$ are in
$\mathbb{L}_1(\mathbb{R})\cap \mathbb{L}_2(\mathbb{R})$ and have the best
smoothness properties as possible, as functions of $z$.
Under reasonable identifiability assumptions, $S_{\theta^0,g}(\theta)$ is minimum if and only
if $\theta=\theta^0$. We propose to estimate $S_{\theta^0,g}(\theta)$ for all
$\theta\in \Theta$ by a quantity depending on the observations
$\Delta_1,\cdots,\Delta_n$, expecting thus that the argument minimum of the
estimator  converges to the argument minimum of
$S_{\theta^0,g}(\theta)$, say $\theta^0$.

We propose a first estimator of $\theta^0$, say $\widehat{\theta}_1$, constructed by minimizing
$S_{n,1}(\theta)$, a consistent estimator of $S_{\theta^0,g}$ 
where $g$ is replaced by a kernel deconvolution estimator. 
We show that under classical assumptions, this estimator is a consistent
estimator of $\theta^0$. Its rate of convergence depends on the
smoothness of $f_\varepsilon$ and on the smoothness on $W(z)f_\beta(z)$, as a
function of $z$. More precisely, its depends on the behavior of the ratios of the Fourier transforms
$(Wf_\beta)^*(t)/\overline{f_\varepsilon^*}(t)$ and
$(Wf_\beta^2)^*(t)/\overline{f_\varepsilon^*}(t)$ as $t$ tends to infinity.
We give upper bounds for the risk of $\widehat{\theta}_1$ for various relative
risks and various types of error
density and  derive sufficient conditions ensuring
the $\sqrt{n}$-consistency and the asymptotic normality. These upper bounds and these
sufficient conditions are deeply studied through examples.
 In particular we show that $\widehat{\theta}_1$ is a $\sqrt{n}$-consistent
asymptotically Gaussian estimator of $\theta^0$ in the Cox model, in the model of excess relative
risk, and when $f_{\beta}$ is a general
polynomial function.

The estimation procedure is related to the problem of the estimation $S_{\theta^0,g}(\theta)$.
 Under conditions ensuring that it can be estimated at the
parametric rate, we propose a second estimator $\widehat{\theta}_2$ which
is $\sqrt{n}$-consistent and asymptotically Gaussian of $\theta^0$.
Clearly, these conditions are not always fullfilled and $\widehat{\theta}_2$
does not always exist, whereas $\widehat{\theta}_1$ can be constructed and studied in all setups.

The paper is organized as follows. Section 2 presents the model and  
the assumptions. In Sections 3 and 4 we present the two estimators and
their asymptotic properties illustrated in Section 5. In Section 6, we comment the
 use of the least squares criterion.
The proofs are gathered in Section 7 and in the Appendix.

\section{Model, assumptions and notations}
\setcounter{equation}{0}
\setcounter{lem}{0}
\setcounter{theo}{0}
Before we describe the estimation procedure, we give notations
used throughout the paper and assumptions commonly done in survival data analysis.

\textbf{Notations}
For two complex-valued functions $u$ and $v$ in
$\mathbb{L}_2(\mathbb{R})\cap \mathbb{L}_1(\mathbb{R})$, let
$$
u^*(x)=\int e^{itx}u(t)dt,  \quad u\star v(x)=\int u(y)v(x-y)dy, \
\text{and} \quad  <u,v>=\int u(x)\overline{v}(x)dx$$ with
$\overline{z}$ the conjugate of a complex number $z$. We also use
the notations $$ \|u\|_1=\int |u(x)|dx, \quad \|u\|^2=\int |u(x)|^2
dx,  \quad\|u\|_\infty=\sup_{x \in
\mathbb{R}}|u(x)|,$$
 and for $\theta \in \mathbb{R}^d$, $$\parallel \theta\parallel_{\ell^2}^2=\sum_{k=1}^d
\theta_k^2.$$
For a map $$\displaystyle\begin{array}{ll}\varphi_\theta~: &\Theta\times \mathbb{R}\longrightarrow \mathbb{R}\\
& (\theta,u) \mapsto \varphi_\theta(u),
\end{array}$$ whenever they exist, the first and
second derivatives with respect to $\theta$ are denoted by
\begin{eqnarray*}
\varphi^{(1)}_{\theta}(\cdot)&=&\left(\varphi_{\theta,j}^{(1)}(\cdot)\right)_j \mbox{
with }\varphi_{\theta,j}^{(1)}(\cdot)=\frac{\partial \varphi_{\theta}(\cdot)}{\partial
  \theta_j}\mbox{ for
}j\in \{1,\cdots,m+p\}\\
\mbox{ and
}\hspace{0.5cm}\varphi^{(2)}_{\theta}(\cdot)&=&\left(\varphi_{\theta,j,k}^{(2)}(\cdot)\right)_{j,k}
\mbox{ with }\varphi_{\theta,j,k}^{(2)}(\cdot)=\frac{\partial^2
  \varphi_{\theta}(\cdot)}{\partial \theta_j\theta_k},\mbox{ for
}j,k\in \{1,\cdots,m+p\}.
\end{eqnarray*}
Throughout the paper $\mathbb{P}$, $\mathbb{E}$ and $\mbox{Var}$ denote
respectively the probability,
the expectation,
and the variance when the underlying
and unknown true parameters are $\theta^0$ and $g$.
Finally we use the notation  $a_-$ for the negative part of $a$, which equals $a$ if
$a\leq 0$ and 0 otherwise.

\textbf{Model assumptions}
\begin{align} 
&\hypothese{rbase}
\mbox{The function } \eta_{\gamma^0} \mbox{ is non-negative and  integrable on
}[0,\tau]
.
\\
&\hypothese{hyp1} 
\mbox{Conditionnally on }Z\mbox{ and }U, \mbox{ the failure time } T \mbox{ and the censoring time } C 
\mbox{ are }\\&\notag \mbox{independent}.
\\
&\hypothese{hyp3}
\mbox{The distribution of the censoring time }C, \mbox{ conditional on }Z \mbox{ and }U, \mbox{ does not }\\ 
& \mbox{ depend on }Z \mbox{ and  } U.\notag
\\
\hypothese{hyp2}
&\mbox{The distribution of the failure time } T, \mbox{ conditional on } Z
\mbox{ and } U,\mbox{ does not   ~~~~}\\ \notag &\mbox{depend on }U.
\end{align}
These assumptions are usual in most frameworks dealing with survival data 
analysis and covariate measured with error, see Andersen \textit{ et al.} \citeyear{ABGK}, Prentice and
Self \citeyear{PRENTICESELF}, 
Prentice\citeyear{PRENTICE82}, Gong \citeyear{GONG90} and Tsiatis \citeyear{TSIATIS}. 
Assumptions \textbf{\eref{hyp1}} and \textbf{\eref{hyp3}} state  that a 
general censorship model is considered, where the censoring time has
an arbitrary distribution independent of the covariates. 
Assumption \textbf{\eref{hyp2}} states that the failure time is independent 
of the observed 
covariate when  the observed and true covariates are both given, i.e. the measurement error is not prognostic. 

We  define the filtration
$$\mathcal{F}_t=\sigma\{Z, U, N(s), \ind_{X \geq s}, 0 \leq s \leq t \leq \tau\}.$$
The intensity of the censored process $N(t)$ with respect to the filtration  
$\mathcal{F}_t$
equals
\begin{eqnarray}\label{intensite}
\lambda(t,\theta^0,Z)=\eta_{\gamma^0}(t)Y(t)f_{\beta^0}(Z).
\end{eqnarray}

It follows from \eref{intensite} 
and from the independence of the observations  
$\Delta_i$, that for the individual $i$ the intensity  and the
compensator process
of the censored process 
$N_i(t)=\ind_{X_i \leq t, D_i=1}$ with respect to the
filtration $\mathcal{F}_t$ are respectively
\begin{eqnarray}\lambda_i(t,\theta^0,Z_i)=\eta_{\gamma^0}(t)Y_i(t)f_{\beta^0}(Z_i)\mbox{
and }\Lambda_i(t,\theta^0,Z_i)=\int_0^t \lambda_i(s,\theta^0,Z_i)ds.
\label{deflambda}
\end{eqnarray}
Moreover the process $M_i(t)=N_i(t)-\Lambda_i(t,\theta^0,Z_i)$ 
is a local square integrable martingale.
As a consequence, the least squares criterion defined in \eref{S0} can be
rewritten as
\begin{multline}
\label{carre}
S_{\theta^0,g}(\theta)=\int_0^\tau
\mathbb{E}\left[\left\lbrace\eta_{\gamma}(t)f_{\beta}(Z)
-\eta_{\gamma^0}(t)f_{\beta^0}(Z)\right\rbrace^2Y(t)W(Z)\right]dt\\
-\int_0^\tau
\mathbb{E}\left[\left\lbrace
\eta_{\gamma^0}(t)f_{\beta^0}(Z)\right\rbrace^2Y(t)W(Z)\right]dt.
\end{multline}
Since we consider general relative risk functions we assume the below minimal smoothness
conditions with respect to $\theta$.

\textbf{Smoothness assumptions}
\begin{align}
&\hypothese{deriv}\mbox{The functions }\beta\mapsto f_\beta
\mbox{ and }\gamma\mapsto \eta_\gamma\mbox{ admit continuous derivatives up to
order 3 }\\&\notag \mbox{with respect to }\beta \mbox{ and }\gamma \mbox{ respectively}.
\end{align}
We denote by $S_{\theta^0,g}^{(1)}(\theta)$ and $S_{\theta^0,g}^{(2)}(\theta)$ the first and
second derivatives of $S_{\theta^0,g}(\theta)$ with respect to
$\theta$. For all $t$ in $[0, \tau]$,  set
$S^{(2)}_{\theta^0,g}(\theta, t)$ the second derivative of
$S_{\theta^0,g}$ 
when the integral is taken over $[0, t]$, with the
convention that 
$S^{(2)}_{\theta^0,g}(\theta)=S^{(2)}_{\theta^0,g}(\theta, \tau)$.

\textbf{Identifiability and moment assumptions}
\begin{align}
\hypothese{maxunique}
&S^{(1)}_{\theta^0,g}(\theta)=0\mbox{ if and only if }\theta=\theta^0. \\
&\mbox{ For all }t\in [0,\tau], \mbox{ the matrix }
S^{(2)}_{\theta^0,g}(\theta^0,t)\mbox{ exists and is positive definite}.
\hypothese{concls}\\
&\hypothese{moment}\mbox{ The quantity }\mathbb{E}(f_{\beta}^2(Z)W(Z))\mbox{
is finite }.
\\
&\hypothese{R2b}\mbox{ For }j=1,\cdots, m,~~\mathbb{E}\vert f_{\beta^0}(Z) f_{\beta^0,j}^{(1)}(Z)W(Z)\vert ^3,~~\mathbb{E}\vert f_{\beta^0}(Z)W(Z)\vert ^3 \mbox{
are finite}.
\end{align}

We denote by $\mathcal{G}$ the set of densities $g$ such that the 
assumptions \eref{hyp1}-\eref{hyp2},\eref{maxunique}-\eref{R2b} hold.

\section{Construction and study of the first  estimator $\widehat{\theta}_1$}
\setcounter{equation}{0}
\setcounter{lem}{0}
\setcounter{theo}{0}
\subsection{Construction}
If the $Z_i$'s were observed, $S_{\theta^0,g}(\theta)$ would be estimated by 
\begin{eqnarray}
\label{Snt}
\quad\tilde S_n(\theta)=-\frac{2}{n}\sum_{i=1}^n
f_{\beta}(Z_i)W(Z_i)\int_0^\tau\!\!\!\eta_{\gamma}(t)dN_i(t)+\frac{1}{n}\sum_{i=1}^nf_{\beta}^2(Z_i)W(Z_i)\int_0^\tau \!\!\!\eta_{\gamma}^2(t)Y_i(t)dt
\end{eqnarray}
and $\theta^0$ would be estimated by minimizing 
$\tilde S_n(\theta)$. Since the $Z_i$'s are unobservable and $Z_i$
independent of $\varepsilon_i$, the density $h$ of $U_i$ equals $h=g\star
f_\varepsilon$. We thus estimate
$S_{\theta^0,g}$ by
\begin{multline}
\label{Sn1}
S_{n,1}(\theta)=-\frac{2}{n}\sum_{i=1}^n (f_{\beta} W)\star
K_{n,C_n}(U_i)\int_0^\tau\!\!\!\eta_{\gamma}(t)dN_i(t)+\frac{1}{n}\sum_{i=1}^n
(f_{\beta}^2 W)\star K_{n,C_n}(U_i)\int_0^\tau \!\!\!\eta_{\gamma}^2(t)Y_i(t)dt,
\end{multline}
where  $K_{n,C_n}(\cdot)=C_n K_n(C_n\cdot)$ is a deconvolution kernel defined
via its Fourier transform, such that $\int {K_n}(x)dx=1$, and
\begin{eqnarray}
\label{kerdeconv}
K_{n,C_n}^*(t)=\frac{K_{C_n}^*(t)}{\overline{f_\varepsilon^*}(t)}=\frac{K^*(t/C_n)}{
\overline{f_\varepsilon^*}(t)},
\end{eqnarray}
with $K^*$
compactly supported satisfying $\vert 1-K^*(t)\vert\leq \ind_{\vert t\vert
\geq 1}$ and $C_n \rightarrow \infty$ as $n\rightarrow \infty$.

Using this criterion we propose to estimate $\theta^0$ by
\begin{eqnarray}
\label{thetac}
\widehat{\theta}_1=\begin{pmatrix}
\widehat{\beta}_1\\\widehat{\gamma}_1\end{pmatrix}=\arg\min_{\theta=(\beta,\gamma)^\top\in\Theta}S_{n,1}(\theta).
\end{eqnarray}
We require for the construction of $S_{n,1}(\theta)$, that
\begin{align}
&\mbox{ the density } f_{\varepsilon} \mbox{ belongs to
  }\mathbb{L}_2(\mathbb{R})\cap \mathbb{L}_\infty(\mathbb{R}) \mbox{ and for all }
x\in \mathbb{R},\,f_{\varepsilon}^*(x)\not=0\hypothese{fepsnn}.
\end{align}

\subsection{Asymptotic properties of the first $\widehat{\theta}_1$}
Assume that
\begin{align}
&\hypothese{R2}\sup_{g\in \mathcal{G}}\parallel f_{\beta^0}^2 g\parallel_2^2\leq
C_2(f_{\beta^0}^2),~~~\sup_{g\in \mathcal{G}}\parallel f_{\beta^0}g\parallel_2^2\leq
C_2(f_{\beta^0}).\\
&\hypothese{cw1} \sup_{\beta\in \mathbb{B}}(W f_{\beta}),
~~ W \mbox{ and } \sup_{\beta\in \mathbb{B}
  } ( W
f^2_{\beta}) \mbox{ belong to }\mathbb{L}_1(\mathbb{R}).\\
&\hypothese{cw11} \sup_{\beta\in \mathbb{B}}(W f^{(1)}_{\beta})\mbox{ and } \sup_{\beta\in \mathbb{B}
  } ( W f_\beta
f^{(1)}_{\beta}) \mbox{ belong to }\mathbb{L}_1(\mathbb{R}).
\end{align}
As in density deconvolution, or for the estimation of the regression function in
errors-in-variables models, the rate of convergence for estimating $\theta^0$ is given by both the
smoothness of $f_\varepsilon$ and the smoothness
of $(f_{\beta}W)(z)$, and $\partial
(f_{\beta} W)(z)/\partial \beta$, as functions of $z$. 
The smoothness of the error density $f_{\varepsilon}$ is described by the decrease of its Fourier transform.
\begin{align}
&\hypothese{condfeps} \mbox{ There exist positive constants
}\underline{C}(f_{\varepsilon}),\overline{C}(f_{\varepsilon}),\mbox{ and
nonnegative }\delta,~
\rho, \alpha \mbox{ and }u_0  \mbox{ such }\\&\notag \mbox{ that }
\underline{C}(f_{\varepsilon}) \leq
\left| f_\varepsilon^*( u) \right| \left| u\right| ^{\alpha} \exp \left( \delta \left|
u\right| ^{\rho}\right) \leq \overline{C}(f_{\varepsilon})\mbox{ for all }\vert
u\vert \geq u_0.
\end{align}
If $\rho=0$, by convention $\delta=0$.
When $\rho=0=\delta$ in (\ref{condfeps}), $f_\varepsilon$ is called "ordinary
smooth".  When $\delta>0$ and $\rho>0$, it is called "super smooth".
Densities satisfying \textbf{\eref{condfeps}} with $\rho>0$ and
$\delta>0$ are infinitely differentiable. The standard examples for
super smooth densities are the Gaussian or Cauchy
distributions which are  super smooth of  respective order $\alpha=0, \rho=2$ and
$\alpha=0, \rho=1$. For ordinary smooth densities,
one can cite for instance the double exponential (also called
Laplace) distribution with $\rho=0=\delta$ and $\alpha=2$. 
We consider here that $0\leq \rho\leq 2$.
The square integrability of $f_{\varepsilon}$ in \eref{fepsnn} requires that
$\alpha> 1/2$ when $\rho=0$ in 
\textbf{\eref{condfeps}}.

The smothness of $f_\beta W$ is described by the following assumption.
\begin{align}
&\hypothese{super} \mbox{ There exist positive constants
}\underline{L}(f),\overline{L}(f) \mbox{ and }a, d, u_0, r \mbox{  nonnegative
numbers such}\\
&\notag\mbox{ that for all }\beta\in \mathbb{B}, f_{\beta}
W \mbox{ and } f_{\beta}^2W\mbox{ and their derivatives up to order 3 with
}
\\ &\notag \mbox{ respect to }\beta, \mbox{ belong to }\\
&\mathcal{H}_{a,d,r}=\{ f \in \mathbb{L}_1(\mathbb{R});
\underline{L}(f)\leq |f^*(u)|\vert u\vert^a\exp( d |u|^r) \leq
\overline{L}(f)<\infty\mbox{  for all }\vert u\vert \geq u_0\}.
\end{align}
If $r=0$, by convention $d=0$.
\begin{theo}
\label{thv1}
Let \textbf{\eref{rbase}}-\textbf{\eref{super}}
hold. Let
 $\widehat{\theta}_1=\widehat{\theta}_1(C_n)$ be defined by \eref{Sn1} and
\eref{thetac} with $C_n$ a sequence such that 
\begin{eqnarray}
\label{condcons2}C_n^{(2\alpha-2a+1-\rho+(1-\rho)_-)}\exp\{-2dC_n^r+2\delta
C_n^\rho\}/n=o(1) \mbox{ as }n \rightarrow +\infty.
\end{eqnarray}
 \textbf{1)} Then
$\mathbb{E}(\parallel\widehat{\theta}_1(C_n)-\theta^0\parallel_{\ell^2}^2)=o(1),$
 as $n\rightarrow \infty$ and $\widehat{\theta}_1(C_n)$ is a consistent
 estimator of $\theta^0.$
\\
\textbf{2)} Moreover, 
$\mathbb{E}(\parallel
\widehat\theta_1-\theta^0\parallel_{\ell^2}^2)=O(\varphi_n^2)$ with
$\varphi_n^2=\|(\varphi_{n,j})\|^2_{\ell^2}$,
$\varphi_{n,j}^2=B_{n,j}^2(\theta^0)+V_{n,j}(\theta^0)/n,$ where
$B_{n,j}^2(\theta^0)=\min \lbrace
B_{n,j}^{[1]}(\theta^0),B_{n,j}^{[2]}(\theta^0)\rbrace$,
$V_{n,j}(\theta^0)=\min \lbrace
V_{n,j}^{[1]}(\theta^0),V_{n,j}^{[2]}(\theta^0)\rbrace$, with
\begin{eqnarray*}
B_{n,j}^{[q]}(\theta^0)&=&
\left\|(f_{\beta^0}^2W)^*(K_{C_n}^*-1)\right\|^2_q+\left\|(f_{\beta^0}W)^*
(K_{C_n}^*-1)\right\|^2_q+\left\|\left(f^{(1)}_{\beta^0,j} W\right)^*(K_{C_n}^*-1)\right\|^2_q\\&&+
\left\|\left(f^{(1)}_{\beta^0,j}f_{\beta^0}W\right)^*  (K_{C_n}^*-1)\right\|^2_q,
\\
V_{n,j}^{[q]}(\theta^0)&=&
\left\|(f_{\beta^0}^2W)^*\frac{K_{C_n}^*}{\overline{f_\varepsilon^*}} \right\|^2_q+\left\|(f_{\beta^0}W)^*
\frac{K_{C_n}^*}{\overline{f_\varepsilon^*}} \right\|^2_q+
\left\|\left(f^{(1)}_{\beta^0,j}
W\right)^*\frac{K_{C_n}^*}{\overline{f_\varepsilon^*}} \right\|^2_q
\\&&+
\left\|\left(f^{(1)}_{\beta^0,j}f_{\beta^0}W\right)^*
\frac{K_{C_n}^*}{\overline{f_\varepsilon^*}} 
\right\|^2_q.\end{eqnarray*}
\textbf{3)} Furthermore, the resulting  rate $\varphi_n^2$ is given in
Table~\ref{rates}
\end{theo}

We point out that the rate for estimating $\beta^0$ depends on the smoothness properties
of $\partial(Wf_\beta)(z)/\partial \beta$ and $\partial(Wf_\beta^2)(z)/\partial \beta$ as a function of $z$, 
whereas, the rate for estimating $\gamma^0$
depends on the smoothness properties of $Wf_\beta(z)$ and  $Wf_\beta^2(z)$ as a
function of $z$. In both cases, the smoothness properties of
$\eta_{\gamma}$  as a function of $t$ does not have influence on the rate of convergence.

The terms $B_{n,j}^2$ and $V_{n,j}$ are respectively the squared bias and variance
terms. As usual ,the bias is the smallest for the  smoothest functions $(W f_\beta)(z)$
and $\partial (f_\beta W)(z)/\partial \beta $, as functions of $z$. As in
density deconvolution, the biggest variance are obtained for the smoothest error
density $f_\varepsilon$. Hence, the slowest rates are
obtained for the smoothest errors density $f_\varepsilon$, for instance for
Gaussian $\varepsilon$'s. 

The rate of convergence of the estimator $\widehat{\theta}_1$ could be improved by assuming smoothness properties on the density $g$.
But, since $g$ is unknown, we choose to not assume such properties.
Consequently, without any additional assumptions on $g$, the parametric rate of convergence is
achieved as soon as $(W
f_\beta)$ and $(W f_\beta^2)$ and their derivatives, as functions
of $z$, are smoother than the
errors density $f_\varepsilon$.

%%%%%%%%%%%%%%%%%%%%%%%%%%%%%%%%%%%%%%%%%%%%%%%%%%%%%
\begin{center}
{\Small
\begin{table}[!ht]
\begin{tabular}{|c|l||c|c||}\cline{3-4}\cline{3-4}
\multicolumn{2}{c||}{} &\multicolumn{2}{c||}{$f_\varepsilon$} \\\cline{3-4}
\multicolumn{2}{c||}{} & $\rho=0$ in \textbf{\eref{condfeps}}& $ \rho>0$ in \textbf{\eref{condfeps}} \\
\multicolumn{2}{c||}{} & ordinary smooth & super smooth \\\hline\hline
\multirow{4}{.7cm}{\\ \vfill\null $Wf_{\beta^0}$} & $\;$ & $\;$ & $\;$ \\
& $\begin{array}{l}
  d=r=0 \\
\mbox{in \textbf{\eref{super}}}\\
  \small{\mbox{Sobolev}}
\end{array}$ & 
$\begin{array}{cc}
&\\
    a < \alpha+1/2 & n^{-\frac{2 a-1}{2\alpha}} \\
&\\
    \hline
&\\
    a \geq \alpha+1/2 & n^{-1} \\
  \end{array}$
  & $ \left(\log n \right)^{-\frac{2 a-1}{\rho}}$
\\$\;$ & $\;$ & $\;$ & $\;$ \\
\cline{2-4}
%$\;$ & $\;$ & $\;$ & $\;$ \\
& $\begin{array}{l}
  r>0\\
\mbox{in \textbf{\eref{super}}}\\
  \mathcal{C}^\infty
\end{array}$ &  $ n^{-1}$  & $
  \begin{array}{cc}
&\\    r < \rho & (\log n)^{A(a,r,\rho)}
    \exp\left\{-2 d\left(\frac{\log n}{2\delta}
\right)^{r/\rho}     \right\} \\
&\\
\hline

r=\rho &
\begin{array}{ll}
&\\
 d < \delta &  {(\log n)^{A(a,r,\rho)+2 \alpha d/(\delta r)}}n^{-d/\delta}
 \\
& \\
d = \delta, \, a < \alpha +1/2 & {(\log n)^{(2 \alpha -2 a +1)/r}}{n^{-1}} \\
& \\
d = \delta, \, a \geq \alpha+1/2 & n^{-1}\\
& \\
d > \delta & n^{-1}\\
&\\
\end{array}
       \\
\hline
&\\
    r > \rho & n^{-1} \\
&\\
  \end{array} $
\\
\hline\hline
\end{tabular}\\~\\
where $A(a, r, \rho)= (-2a +1 -r +(1-r)_-)/\rho$.
\caption{Rates of convergence $\varphi_n^2$ of $\widehat\theta_1$ }
\label{rates}
\end{table}}
\end{center}
\vspace{-1cm}

\subsection{Consequence : a sufficient condition to obtain the parametric rate
of convergence with $\widehat{\theta}_1$}

\begin{align}
&\condition{C1}\mbox{ There exists a weight function }W
\mbox{ such
  that  the functions }\\&\notag
\sup_{\beta\in \mathbb{B}}(f_{\beta}
W)^*/\overline{f_\varepsilon^*}\,,\,\sup_{\beta\in \mathbb{B}}(f_{\beta}^2
W)^*/\overline{f_\varepsilon^*}   \mbox{ belong to }
\mathbb{L}_1(\mathbb{R})\cap \mathbb{L}_2(\mathbb{R}).\\
&\condition{C2} \mbox{ The functions }\sup_{\beta\in
\mathbb{B}}\Big(f^{(1)}_{\beta} W\Big)^*/\overline{f_\varepsilon^*} \mbox{ and }  \sup_{\beta\in
\mathbb{B}}\Big(f^{(1)}_{\beta}f_{\beta} W 
\Big)^*/\overline{f_\varepsilon^*} \\&\notag \mbox{ belong to }
\mathbb{L}_1(\mathbb{R})\cap \mathbb{L}_2(\mathbb{R}).\end{align}
\begin{align}
&\condition{C3} \mbox{ The functions }\Big(f^{(2)}_{\beta} W 
\Big)^*/\overline{f_\varepsilon^*}  \mbox{ and }  \Big(\frac{\partial^2(f_{\beta} ^2 W)}{\partial\beta^2} 
\Big)^*/\overline{f_\varepsilon^*} \\&\notag \mbox{ belong to }
\mathbb{L}_1(\mathbb{R})\cap \mathbb{L}_2(\mathbb{R})\mbox{ for all }\beta\in
\mathbb{B}.
\end{align}

\begin{theo}
\label{thCS} Let
\textbf{\eref{rbase}}-\textbf{\eref{cw1}} and
\textbf{\eref{C1}}-\textbf{\eref{C3}} hold.
 Then $\widehat{\theta}_1$ defined by \eref{thetac} is a
$\sqrt{n}$-consistent estimator of $\theta^0$. Moreover
$$\sqrt{n}(\widehat{\theta}_1-\theta^0)\cvl \mathcal{N}(0,\Sigma_1),$$ where
$\Sigma_1$ equals 
\begin{multline}
\left(\mathbb{E}\left[
-2\int_0^\tau\frac{\partial^2 ((f_\beta W)(Z)
\eta_{\gamma}(s))}{\partial
\theta^2}\left.\right|_{\theta=\theta^0}dN(s)+
\int_0^\tau \frac{\partial^2 ((f_\beta^2W)(Z) 
\eta^2_{\gamma}(s))}{\partial
\theta^2}\left.\right|_{\theta=\theta^0}Y(s)ds\right]\right)^{-1}\\
\quad\times
\Sigma_{0,1}\left(\mathbb{E}\left[
-2\int_0^\tau \frac{\partial^2 ((f_\beta
W)(Z)
\eta_{\gamma}(s))}{\partial
\theta^2}\left.\right|_{\theta=\theta^0}dN(s)+
\int_0^\tau \frac{\partial^2 ((f_\beta^2
W)(Z)
\eta^2_{\gamma}(s))}{\partial
\theta^2}\left.\right|_{\theta=\theta^0}Y(s)ds\right]\right)^{-1}
\end{multline}
with
\begin{multline*}
\displaystyle
\Sigma_{0,1}=\mathbb{E}\left\lbrace\left[
-2\int_0^\tau \frac{\partial (R_{\beta,f_\varepsilon,1}(U)\eta_{\gamma}(s))}{\partial
\theta}\left.\right|_{\theta=\theta^0}dN(s)+
\int_0^\tau \frac{\partial (R_{\beta,f_\varepsilon,2}(U) \eta^2_{\gamma}(s))}{\partial
\theta}\left.\right|_{\theta=\theta^0}Y(s)ds\right]\right. \\\left.\times
\left[
-2\int_0^\tau \frac{\partial (R_{\beta,f_\varepsilon,1}(U)\eta_{\gamma}(s))}{\partial
\theta}\left.\right|_{\theta=\theta^0}dN(s)+
\int_0^\tau \frac{\partial (R_{\beta,f_\varepsilon,2}(U) \eta^2_{\gamma}(s))}{\partial
\theta}\left.\right|_{\theta=\theta^0}Y(s)ds\right]^\top    \right\rbrace
\end{multline*}
where
$$R_{\beta,f_\varepsilon, 1}(U)=\int
(Wf_\beta)^*(t)\frac{e^{-itU}}{\overline{f_\varepsilon^*}(t)}dt\quad \mbox{and }\quad R_{\beta,f_\varepsilon,2}(U)=\int
(Wf_\beta^2)^*(t)\frac{e^{-itU}}{\overline{f_\varepsilon^*}(t)}dt.$$
\end{theo}

The conditions \textbf{\eref{C1}}-\textbf{\eref{C3}}, stronger than
\textbf{\eref{condfeps}} and \textbf{\eref{super}}, ensure the existence of the functions
$R_{\beta,f_\varepsilon,j}$ for $j=1,2.$

\section{Construction and study of the second estimator $\widehat{\theta}_2$}
\setcounter{equation}{0}
\setcounter{lem}{0}
\setcounter{theo}{0}
\subsection{Construction}
Our estimation procedure, based on the estimation of the least squares
criterion, requires the estimation of  $\mathbb{E}[\int_0^\tau (f_\beta
W)(Z)dN(t)]$ and $\mathbb{E}[\int_0^\tau (f^2_\beta
W)(Z)Y(t)dt]$, which are linear functional of $g$. It may appear that these
linear functional could be directly estimated, without kernel deconvolution plugg-in.
In this context, we propose another estimator of $\theta^0$. It is based on sufficient conditions allowing to construct
 a $\sqrt{n}$-consistent estimator of these linear functionals and hence to
estimate $S_{\theta^0,g}$ with the parametric rate.

We say that the conditions \eref{C21}-\eref{C23} hold if there exist a weight
function $W$ and two functions $\Phi_{\beta,f_\varepsilon,1}$ and $\Phi_{\beta,f_\varepsilon,2}$
not depending on $g$, such that for all $\beta \in \mathbb{B}$ and for all $g$
\begin{align}
&\condition{C21}\, \mathbb{E}_{\theta^0,g}\left[\int_0^\tau (f_\beta
W)(Z)dN(t)\right]= \mathbb{E}_{\theta^0,h}\left[\int_0^\tau \Phi_{\beta,f_\varepsilon,1}(U)dN(t)\right]\\
\notag & \mbox{ and } 
\mathbb{E}_{\theta^0,g}\left[\int_0^\tau (f_\beta^2
W)(Z)Y(t)dt\right]= \mathbb{E}_{\theta^0,h}\left[\int_0^\tau \Phi_{\beta,f_\varepsilon,2}(U)Y(t)dt\right];
\end{align}
\begin{align}
&\condition{C22} \mbox{ For }k=0,1,2\mbox{ and for }j=1,2,\quad\mathbb{E}[\sup_{\beta\in \mathbb{B}}\|\Phi^{(k)}_{\beta,f_\varepsilon,j}(U)\|_{\ell^2}]<\infty;
\end{align}
\begin{align}
&\condition{C23}\mbox{ For }j=1,2\mbox{ and for all }\beta\in \mathbb{B},\,\mathbb{E}\left[
\|\Phi^{(1)}_{\beta,f_\varepsilon,j}(U)\|_{\ell^2}^2\right]
<\infty.
\end{align}
Under \eref{C21}-\eref{C23}, we estimate $S_{\theta^0g}$ by
\begin{eqnarray}
\label{Sn2}
\quad S_{n,2}(\theta)=-\frac{2}{n}\sum_{i=1}^n
\int_0^\tau \Phi_{\beta,f\varepsilon,1}(U_i)\eta_{\gamma}(t)dN_i(t)+\frac{1}{n}\sum_{i=1}^n\int_0^\tau \Phi_{\beta,f\varepsilon,2}(U_i)\eta_{\gamma}^2(t)Y_i(t)dt
\end{eqnarray}
and $\theta^0$ is estimated by
\begin{eqnarray}
\label{thetac2}
\widehat{\theta}_2=\arg\min_{\theta\in \Theta}S_{n,2}(\theta).
\end{eqnarray}

The main difficulty for finding such functions $\Phi_{\beta,f_\varepsilon,1}$
and $\Phi_{\beta,f_\varepsilon,2}$
lies in the constraint that they must not depend on the unknown
density $g$. 
We refer to Section \ref{commentaire} for details
on how to construct such functions $\Phi_{\beta,f_\varepsilon,j}$, $j=1,2$.

\subsection{Asymptotic properties of $\widehat{\theta}_2$}

\begin{theo}
\label{thv2} Let \eref{rbase}-\eref{R2b}, and the conditions
\eref{C21}-\eref{C23} hold.
Then $\widehat{\theta}_2$, defined by \eref{thetac2} is a $\sqrt{n}$-consistent
estimator of $\theta^0$. Moreover $$\sqrt{n}(\widehat{\theta}_2-\theta^0)\cvl
\mathcal{N}(0,\Sigma_2),$$ where $\Sigma_2$ equals 
\begin{multline}
\left(\mathbb{E}\left[
-2\int_0^\tau \frac{\partial^2(\Phi_{\beta,f_\varepsilon,1}(U)\eta_{\gamma}(s))}{\partial
\theta^2}\left.\right\vert_{\theta=\theta^0}dN(s)+
\int_0^\tau \frac{\partial^2 (\Phi_{\beta,f_\varepsilon,2}(U) \eta^2_{\gamma}(s))}{\partial
\theta^2}\left.\right\vert_{\theta=\theta^0}Y(s)ds\right]\right)^{-1}\\
\quad\times
\Sigma_{0,2}\left(\mathbb{E}\left[
-2\int_0^\tau \frac{\partial^2 (\Phi_{\beta,f_\varepsilon,1}(U)\eta_{\gamma}(s))}{\partial
\theta^2}\left.\right\vert_{\theta=\theta^0}dN(s)+
\int_0^\tau \frac{\partial^2 (\Phi_{\beta,f_\varepsilon,2}(U) \eta^2_{\gamma}(s))}{\partial
\theta^2}\left.\right\vert_{\theta=\theta^0}Y(s)ds\right]\right)^{-1}
\end{multline}
with
\begin{multline*}
\displaystyle
\Sigma_{0,2}=\mathbb{E}\left\lbrace\left[
-2\int_0^\tau \frac{\partial (\Phi_{\beta,f_\varepsilon,1}(U)\eta_{\gamma}(s))}{\partial
\theta}\left.\right\vert_{\theta=\theta^0}dN(s)+
\int_0^\tau \frac{\partial (\Phi_{\beta,f_\varepsilon,2}(U) \eta^2_{\gamma}(s))}{\partial
\theta}\left.\right\vert_{\theta=\theta^0}Y(s)ds\right]\right. \\\left.\times
\left[
-2\int_0^\tau \frac{\partial (\Phi_{\beta,f_\varepsilon,1}(U)\eta_{\gamma}(s))}{\partial
\theta}\left.\right\vert_{\theta=\theta^0}dN(s)+
\int_0^\tau \frac{\partial (\Phi_{\beta,f_\varepsilon,2}(U) \eta^2_{\gamma}(s))}{\partial
\theta}\left.\right\vert_{\theta=\theta^0}Y(s)ds\right]^\top    \right\rbrace.
\end{multline*}
\end{theo}

\subsection{Comments on conditions ensuring $\sqrt{n}$-consistency~:
comparison of $\widehat{\theta_1}$ and $\widehat{\theta}_2$}

\label{commentaire}
{\rm
~~\\

Let us briefly compare the conditions \eref{C1}-\eref{C3} to the
conditions \eref{C21}-\eref{C23}.
It  is  noteworthy  that   the  conditions  \eref{C21}-\eref{C23}  are  more
general. First, the condition \eref{C21} does not require that $f_\beta W$,
$f_\beta^2 W$ belong to
$\mathbb{L}_1(\mathbb{R})$ (as for instance in the Cox Model).
Second, we point out that Condition \eref{C1} implies \eref{C21},
with $ \Phi_{\beta,f_\varepsilon,j}=R_{\beta,f_\varepsilon,j}$. 
This comes from the facts that under \eref{C1}-\eref{C3}, by denoting $ \Phi_{\beta,f_\varepsilon,1}=
(W f_\beta)^*/\overline{f_\varepsilon^*}$ and $\Phi_{\beta,f_\varepsilon,2}
=(W f^2_\beta)^*/\overline{f_\varepsilon^*}$, we have
\begin{eqnarray*}
\mathbb{E}[Y(t)\Phi_{\beta,f_\varepsilon,2}(U)]
&=&\iiint \ind_{x\geq
t}\Phi_{\beta,f_\varepsilon,2}(u)f_{X,Z}(x,z)f_\varepsilon(u-z)dx\,du\,dz\\
&=&
\iint \ind_{x\geq
t}f_{X,Z}(x,z)\frac{1}{2\pi}\int
\Phi^*_{\beta,f_\varepsilon,2}(s)e^{-isz}\overline{f_\varepsilon^*}(s)ds\,dx\,dz\\&=&
\iint \ind_{x\geq
t}f_{X,Z}(x,z)\frac{1}{2\pi}\int
\frac{(W
f_\beta^2)^*(s)}{\overline{f_\varepsilon^*}(s)}e^{-isz}\overline{f_\varepsilon^*}(s)ds\,dx\,dz\\
&=&\mathbb{E}[Y(t)(W f_\beta^2 )(Z)].\end{eqnarray*}
Consequently
\begin{eqnarray*}
\mathbb{E}\Big[\int_0^\tau\Phi_{\beta,f_\varepsilon,2}(U)\eta_{\gamma}^2(t)Y(t)dt\Big]
&=&\mathbb{E}\Big[\int_0^\tau f_{\beta}^2(Z)W(Z)\eta_{\gamma}^2(t)dt\Big],
\end{eqnarray*}
 and analogoulsy
\begin{eqnarray*}
\mathbb{E}\Big[\int_0^\tau\Phi_{\beta,f_\varepsilon,1}(U)\eta_{\gamma}(t)dN(t)\Big]
&=&\mathbb{E}\Big[\int_0^\tau f_{\beta}(Z)W(Z)\eta_{\gamma}(t)dN(t)\Big].
\end{eqnarray*}
Hence Condition \eref{C21} holds and
$\Sigma_{0,1}=\Sigma_{0,2}$ with $\Sigma_{0,1}$ defined in
Theorem  \ref{thCS}.

These comments underline the key importance of the weight function $W$.
For instance, if $f_\beta(z)=1-\beta+\beta/(1+z^2)$, and
$f_\varepsilon$ is the Gaussian density, then 
it seems impossible to find a function $\Phi_{\beta,f_\varepsilon,2}$ such
that
$\mathbb{E}[Y(t)\Phi_{\beta,f_\varepsilon,2}(U)]=\mathbb{E}[Y(t)f_\beta^2(Z)],$
whereas \eref{C1}-\eref{C3} hold 
by taking
$W(z)=(1+z^2)^4\exp(- z^2/(4 \delta))$. 
In this special example, we exhibit a  suitable choice of $W$ that ensures
that
 condition \eref{C1}-\eref{C3} are fulfilled (see Section~\ref{exemples} for further details).
Nevertheless, such weight function are not always available and hence those
conditions \eref{C1}-\eref{C3} are not always fulfilled.
}

\section{Examples}
\label{exemples}
\setcounter{equation}{0}
\setcounter{lem}{0}
\setcounter{theo}{0}

In this section, we illustrate the asymptotic
properties of $\widehat{\theta}_1$ and $\widehat{\theta}_2$ for various
relative risks and error density
$f_\varepsilon$. In all of
these examples, $K^*(t)= \ind_{\vert t\vert \leq 1}$ and the noise distribution is arbitrary, as far as it satisfies
\textbf{\eref{fepsnn}} and \textbf{\eref{condfeps}} with $0\leq\rho\leq 2$.

The first example deals with Cox model. We show that our
estimation procedure,
based on a nonparametric method and specifically on density deconvolution, also
provides  $\sqrt{n}$-consistent and asymptotically Gaussian
estimator of $\beta^0$. 
The aim of this example is to show that we recover
previous known results using estimators that are quite different from the ones  
proposed by Nakamura \citeyear{NAKAMURA92} and studied by Kong and Gu
\citeyear{KONGGU} or from the ones  proposed by Augustin
\citeyear{Augustin04}.

The others examples we consider, deal with
relative risks for which no consistent estimators were known when the covariate
is mismeasured.

\begin{exem}{\rm \textbf{Exponential relative risk (Cox model)}
\label{lemexp1}

Let $f_\beta$ be of the form $f_\beta(z)=\exp (\beta z) $ and assume that
$\mathbb{E}[\exp (2\beta^0 U)]<\infty.$
Let
 $W(z)=\exp\{-z^2/(4 \delta)\}$. Then the conditions \eref{C1}-\eref{C3}
as well as the condition \eref{C21}-\eref{C23} are satisfied. Consequently
the estimators $\widehat{\theta}_1$ and $\widehat{\theta}_2$ are
$\sqrt{n}$-consistent and asymptotically Gaussian estimators of $\theta^0$, with the same asymptotic
variance.

One could also choose $W\equiv 1$ and use that
$\mathbb{E}[\exp (\beta Z )]=\mathbb{E}[ \exp ( \beta
U)]/\mathbb{E}[\exp(\beta\varepsilon)] .$
This implies that if we denote by
\begin{eqnarray*}
\Phi_{\beta,f_\varepsilon,1}(U)
=\frac{ \exp(2\beta U )}{\mathbb{E}[\exp(2\beta \varepsilon)]}
\mbox{ and }\Phi_{\beta,f_\varepsilon,2}(U)=\frac{
\exp(\beta U )}{\mathbb{E}[\exp(\beta \varepsilon)]}
\end{eqnarray*}
then $
\mathbb{E}[\Phi_{\beta,f_\varepsilon,1}(U)]
=\mathbb{E}[f_\beta^2(Z)]$ and $\mathbb{E}[Y(t)\Phi_{\beta,f_\varepsilon,2}(U)]
=\mathbb{E}[Y(t) f_\beta(Z)], $
and the criterion $S_{n,2}$ defined by \eref{Sn2} exists.

In this case $\widehat{\theta}_2$ is also a $ \sqrt{n}$-consistent and
asymptotically Gaussian estimator of $\theta^0$. 
}
\end{exem}

\begin{exem}{\bf Polynomial relative risk 1 (included Excess relative risk model)}
\label{lempol1}
{\rm
Let $f_\beta$ be of the form $f_\beta(z)=1+\sum_{k=1}^m\beta_k z^k$ and let
$W(z)=\exp\{-z^2/(4 \delta)\}$. Then conditions \eref{C1}-\eref{C3}
as well as conditions \eref{C21}-\eref{C23} are satisfied. Consequently
the estimators $\widehat{\theta}_1$ and $\widehat{\theta}_2$ are
$\sqrt{n}$-consistent and asymptotically Gaussian estimators of $\theta^0$,
with the same asymptotic variance.

We point out that when $m=1$, $f_\beta(z)=1+\beta z$, and this model is known as the model of excess relative risk.

One can also choose $W\equiv 1$, provided that
  the kernel $K$ has finite absolute moments of order $m$ and satisfies
  $\int u^r K(u)du=0$, for $r=1,\ldots,m$.
With this choice of $W$, $\widehat{\theta}_1$  remains a $ \sqrt{n}$-consistent and
asymptotically Gaussian estimator of $\theta^0$. }
\end{exem}
\begin{exem}{\bf Cosines relative risk 1}
\label{lemcos}
{\rm
Let $f_\beta$ be of the form $f_\beta(z)=\sum_{j=1}^m \beta_j\cos(jz)$ with
$\sum_{j=1}^m\beta_j=1$.
Let
 $W(z)=\exp\{-z^2/(4 \delta)\}$. Then the conditions \eref{C1}-\eref{C3}
 as well as conditions \eref{C21}-\eref{C23} are satisfied. Consequently
the estimators $\widehat{\theta}_1$ and $\widehat{\theta}_2$ are
$\sqrt{n}$-consistent and asymptotically Gaussian estimators of $\theta^0$, with the same asymptotic
variance.

One can also choose $W\equiv 1$ and use that
$\mathbb{E}[\exp (ij Z)] = \mathbb{E}[\exp (ijU)]/\mathbb{E}[\exp(ij\varepsilon)].
$
This implies that if we denote  by
\begin{eqnarray*}
\Phi_{\beta,f_\varepsilon,1}(U)=\frac{1}{2}\left[\frac{
\exp(ijU  )}{f_\varepsilon ^*(j) }+\frac{\exp(-ijU  )}{\overline{f_\varepsilon^*}(j)}\right]
\end{eqnarray*}
and 
\begin{multline*}
\Phi_{\beta,f_\varepsilon,2}(U)
=\frac{1}{4}\left\lbrace 1+\sum_{j=1}^m\beta_j^2
\left[\frac{ \exp(2ij U )}{f_\varepsilon^*(2j)}+\frac{\exp(-2ijU)}{\overline{f_\varepsilon^*}(2j)}\right]
\right.\\
+\left.\sum_{j=1}^m\sum_{k\not=j}\beta_j\beta_k\left[\frac{\exp(i(j+k)U)}{f_\varepsilon^*(j+k)}
+\frac{\exp(-i(j+k)U)}{\overline{f_\varepsilon^*}(j+k)}+\frac{\exp(i(j-k)U)}{f_\varepsilon^*(j-k)}
+\frac{\exp(i(-j+k)U)}{\overline{f_\varepsilon^*}(j-k)}\right]\right\rbrace
\end{multline*}
 then the criterion $S_{n,2}$ defined in \eref{Sn2} exists.

With this choice of $W$, $\widehat{\theta}_2$ remains a $ \sqrt{n}$-consistent and
asymptotically Gaussian estimator of $\theta^0$. In the same way,
$\widehat{\theta}_1$  with $W\equiv
1$ also remains $ \sqrt{n}$-consistent and
asymptotically Gaussian estimator of $\theta^0$.  }
\end{exem}
\begin{exem}{\bf Cauchy relative risk 1}
\label{cauchy1}
{\rm
Consider $f_\beta$ of the form $f_\beta(z)=1-\beta+\beta/(1+z^2)$. Then
$f_\beta$ has the regularity  of $z\mapsto 1/(1+z^2)$ which belongs to 
$\mathcal{H}_{a,d,r}$ defined in \textbf{\eref{super}}
with $a=0, \, d=1/2$ and $r=1$.  Let $W(z)=(1+z^2)^4 \exp\{-z^2/(4 \delta)\}$.
Hence
the functions $f_\beta W$,
$f_\beta^2 W$ and their derivatives in $\beta$ up to order 3 belong to
$\mathcal{H}_{a,d,r}$ defined in \textbf{\eref{super}}
with $\rho <r=2$ or $\rho=r=2$ and $d>\delta$. Consequently, the conditions \eref{C1}-\eref{C3}
 as well as conditions \eref{C21}-\eref{C23} are satisfied and
the estimators $\widehat{\theta}_1$ and $\widehat{\theta}_2$ are
$\sqrt{n}$-consistent and asymptotically Gaussian estimators of $\theta^0$, with the same asymptotic
variance.}
\end{exem}

{\rm This simple example underlines the importance of the smoothing weight
function $W$ in the construction of $\widehat{\theta}_1$ or
$\widehat{\theta}_2$. Indeed, without a smoothing function $W$ in front of the relative risk,
Theorem~\ref{thv1} predicts a rate
of convergence of order $\exp(-2\sqrt{\log n})$ for Gaussian $\varepsilon$. }

 \begin{exem}{\bf Laplace relative risk}
\label{symexp}
{\rm
Consider $f_\beta$ of the form $f_\beta(z)= 1+\beta f(z)$ with $f(z)= \exp
(-|z|/2)-1$. Since the Fourier transform of $z \mapsto \exp
(-|z|/2)$ is slowly decaying, like $|u|^{-2}$ as $|u| \to \infty$, if we
choose $W \equiv 1$, the estimator $\widehat \theta_1$
is not $\sqrt{n}$-consistent  as soon as $|f_\varepsilon ^*(u)|\leq o(|u|^{-2})$ with $|u| \to \infty$. A closer look tells us
that $f_\beta$ and its derivative in $\beta$ is $\mathcal{C}^\infty$ except at one point $z=0$.
Therefore, a proper choice of $W$ can smooth out at 0 and make $W f_\beta$, $W f_\beta^2$ and their
derivatives in $\beta$ infinitely differentiable functions in $z$. This choice
of $W$ ensures the $\sqrt{n}$-consistency of $\widehat \theta_1$ whatever
$f_\varepsilon$ satisfies  \textbf{\eref{condfeps}} with $0 < \rho < 1$. Even
if $\rho\geq 1$, the rate of
$\widehat \theta_1$ is much faster when using our choice of $W$ then it would be
for $W \equiv 1$. Let us precise the choice of $W$. Set \begin{equation}\label{smoothphi}
\Psi_{A,B,R}(z) = \exp\left(-\frac{1}{(z-A)^R (B-z)^R} \right) I_{[A,B]}(z),
\end{equation}
where $-\infty < A < B < \infty$ are fixed and $ R>0$. According to
Lepski and Levit~\citeyear{LepskiLevit} and Fedoryuk~\citeyear{Fedoryuk},
p. 346, Theorem 7.3, $
|\Psi^*_{A,B,R}(u)| \leq c\exp(-C |u|^{R/(R+1)}), \text{ as } |u|\to \infty
$
and $c,\,C$ are positive constants.
We propose to take $W$ equal to $\Psi_{0,100,R}$ or
$\Psi_{-100,0,R}$ or their sum.

This choice of $W$ ensures that $f_\beta W$, $f_\beta^2 W$ and their derivatives up to order 3
belong to $\mathcal{H}_{a,d,r}$ defined in \textbf{\eref{super}} with  $d>0$
and 
$r=R/(R+1)$ closer to 1 as $R$ comes larger.

If $f_\varepsilon$ satisfies \textbf{\eref{condfeps}} with $0 \leq \rho < 1$,
we choose $R$ large enough
such that $r=R/(R+1) > \rho$. Hence, the conditions \eref{C1}-\eref{C3}
as well as the conditions \eref{C21}-\eref{C23} are satisfied. Consequently
the estimators $\widehat{\theta}_1$ and $\widehat{\theta}_2$ are
$\sqrt{n}$-consistent and asymptotically Gaussian estimators of $\theta^0$, with the same asymptotic
variance.

If $\rho\geq 1$, for this choice of $W$, the functions $Wf_\beta$ and  $Wf_\beta^2$ and their derivatives in $\beta$ up to order 3,
belong to $\mathcal{H}_{a,d,r}$ with $r = R/(R+1)$ and hence, according to Table~\ref{rates}, 
$$\mathbb{E}\parallel\widehat{\theta}_{1} -\theta^0\parallel_{\ell^2}^2
=O(1)\left(\log n \right)^{\frac{1-2a-r}{\rho}}
\exp\{-2d (\log n/(2\delta))^{r/\rho}\}.$$
}
\end{exem}

\begin{exem}{\bf Irregular relative risk}
\label{uniform}
{\rm
Consider $f_\beta$ of the form $f_\beta(z)=1-\beta +\beta\ind_{[-1,1]}(z)$ and
take $W=\Psi_{-1,1,R}$ defined by
(\ref{smoothphi}) for $R>0$ .

If $\rho=0$ in \textbf{\eref{condfeps}}, then Conditions \eref{C1}-\eref{C3}
 as well as Conditions \eref{C21}-\eref{C23} are satisfied. Consequently
the estimators $\widehat{\theta}_1$ and $\widehat{\theta}_2$ are
$\sqrt{n}$-consistent and asymptotically Gaussian estimators of $\theta^0$, with the same asymptotic
variance. 

If $\rho>0$, then the best rate for estimating $\theta^0$ is obtained
by choosing $W=\Psi_{-1,1,R}$ with $R>0$ sufficiently large such that
$Wf_\beta$ and  $Wf_\beta^2$ and their derivatives in $\beta$ up to order 3,
belong to $\mathcal{H}_{a,d,r}$ defined in \textbf{\eref{super}} with $0<r  = R/(R+1)<1$ as
close to 1 as needed. 

It follows that if $0\leq\rho<1$, then we can find $W=\Psi_{-1,1,R}$ belonging
to $\mathcal{H}_{a,d,r}$ with $r = R/(R+1)>\rho$. Hence the conditions \eref{C1}-\eref{C3}
 as well as conditions \eref{C21}-\eref{C23} are satisfied and  the
 estimators $\widehat{\theta}_1$ and $\widehat{\theta}_2$ are
$\sqrt{n}$-consistent and asymptotically Gaussian estimators of $\theta^0$.

If $\rho\geq1$, for $W=\Psi_{-1,1,R}$, the functions $Wf_\beta$ and  $Wf_\beta^2$ and their derivatives in $\beta$ up to order 3,
belong to $\mathcal{H}_{a,d,r}$ with $r = R/(R+1)$ and hence, according to
Table~\ref{rates},
\begin{eqnarray*}
\mathbb{E}\parallel\widehat{\theta}_{1} -\theta^0\parallel_{\ell^2}^2
=O(1)\left(\log n \right)^{\frac{1-2a-r}{\rho}}
\exp\{-2d(\log n/(2\delta))^{r/\rho}\}.\end{eqnarray*}
}
\end{exem}

 \begin{exem}{\bf Polygonal relative risk}
\label{polyg}
{\rm
Consider $f_\beta$ with
$f_\beta(z)=1-\beta_2a_--\beta_3|b|^3+\beta_1z+\beta_2(z-a)\ind_{z\geq a}
+\beta_3|z-b|^3$. This relative risk is $\mathcal{C}^\infty$
except at points $a$ and $b$ where it is not differentiable. We suggest to use the
smoothing weight function in (\ref{smoothphi}) as follows. For $R>0$, let
$$
W(z)= \Psi_{a-100,a,R}(z) + \Psi_{a,b,R}(z) + \Psi_{b, b+100,R}(z).
$$

If the noise satisfies \textbf{\eref{condfeps}} with $0 \leq \rho < 1$, then take $R$ large enough
such that $r=R/(R+1) > \rho$ and thus conditions \eref{C1}-\eref{C3}
as well as conditions \eref{C21}-\eref{C23} are satisfied. Consequently
the estimators $\widehat{\theta}_1$ and $\widehat{\theta}_2$ are
$\sqrt{n}$-consistent and asymptotically Gaussian estimators of $\theta^0$, with the same asymptotic
variance.

If $\rho\geq 1$ in \textbf{\eref{condfeps}}, the functions $Wf_\beta$ and  $Wf_\beta^2$ and their derivatives in $\beta$ up to order 3,
belong to $\mathcal{H}_{a,d,r}$ with $r = R/(R+1)$ and hence, according to Table~\ref{rates}
\begin{eqnarray*}
\mathbb{E}\parallel\widehat{\theta}_{1} -\theta^0\parallel_{\ell^2}^2
=O(1)\left(\log n \right)^{\frac{1-2a-r}{\rho}}
\exp\{-2d(\log n/(2\delta))^{r/\rho}\}.
\end{eqnarray*}
}
\end{exem}

\textbf{Comments on the examples \ref{symexp}, \ref{uniform}, and
\ref{polyg}}

{\rm 
In these three examples, $f_\beta W$ belongs to $\mathcal{H}_{a,d,r}$ defined in \textbf{\eref{super}}
with $r$ at most such that $r<1$. Hence $\widehat{\theta}_1$
achieves the $\sqrt{n}$-rate of convergence provided that $f_\varepsilon$ is
ordinary smooth or super smooth with an exponent $\rho<1$.
It
seems therefore impossible to have $(W f_\beta)^*/f_\varepsilon^*$ in
$\mathbb{L}_1(\mathbb {R})$ when the $\varepsilon_i$'s are Gaussian. This comes
from the fact that for these relative risks, the least squares criterion $S_{\theta^0,g}(\theta)$ cannot be estimated with the parametric rate
of convergence and hence could probably, not provide
a $\sqrt{n}$-consistent estimator of $\theta^0$. Nevertheless,
even in cases where the parametric rate of convergence seems not achievable by
such estimators, the resulting rate of the risk of
$\widehat{\theta}_1$ is clearly infinitely faster than the logarithmic rate
predicted by Table \ref{rates} that we
could have without $W$ .}
\bigskip

In most of previous examples where the weight function $W$ is required, the points
where $f_\beta(z)$ has to be smoothed do not depend on $\beta$.
But in survival data analysis the relative risks $f_\beta$ are
usually of the form $f_\beta(z)=f(\beta z)$ (see for instance Prentice and
Self \citeyear{PRENTICESELF}). In such models, the points where
$f_\beta(z)$ has to be smoothed (as function of $z$) will depend on $\beta$.

Let us present such examples.

\begin{exem}{\bf Polynomial relative risk 2}
\label{lempol2}
{\rm
Let $f_\beta$ be of the form $f(\beta z)$ with $f(z)=1+\sum_{k=1}^m a_k z^k$ 
with known $a_k$'s. Let
$W(z)=\exp\{-z^2/(4 \delta)\}$. Then conditions \eref{C1}-\eref{C3}
as well as conditions \eref{C21}-\eref{C23} are satisfied. Consequently
the estimators $\widehat{\theta}_1$ and $\widehat{\theta}_2$ are
$\sqrt{n}$-consistent and asymptotically Gaussian estimators of $\theta^0$,
with the same asymptotic variance.}
\end{exem}

\begin{exem}{\bf Cosines relative risk 2}
\label{lemcos2}
{\rm
Let $f_\beta$ be of the form $f(\beta z)$ with $f(z)=\sum_{j=1}^m a_j\cos(jz)$
with known $a_k$'s such that $\sum_{j=1}^m a_j =1$.
Let
 $W(z)=\exp\{-z^2/(4 \delta)\}$. Then the conditions \eref{C1}-\eref{C3}
 as well as conditions \eref{C21}-\eref{C23} are satisfied. Consequently
the estimators $\widehat{\theta}_1$ and $\widehat{\theta}_2$ are
$\sqrt{n}$-consistent and asymptotically Gaussian estimators of $\theta^0$, with the same asymptotic
variance.}
\end{exem}

\begin{exem}{\bf Cauchy relative risk 2}
\label{cauchy2}
{\rm
Consider $f_\beta$ of the form $f(\beta z)$ with $f(z)=1/(1+z^2)$.
Let
$W(z)=(1+z^2)^4 \exp\{-z^2/(4 \delta)\}$ or $W\equiv 1$. With these choices of $W$,
the functions $f_\beta W$,
$f_\beta^2 W$ and their derivatives in $\beta$ up to order 3 belong to
$\mathcal{H}_{a,d,r}$ defined in \textbf{\eref{super}}
 with $a=0, \, d=1/\beta$ and $r=1$. 
According to Table \ref{rates}, if $f_\varepsilon$ satisfies
\textbf{\eref{condfeps}} with $0\leq \rho<1$, then $\widehat{\theta}_1$ and $\widehat{\theta}_2$ are
$\sqrt{n}$-consistent and asymptotically Gaussian estimators of $\theta^0$, with the same asymptotic
variance. If $f_\varepsilon$ satisfies \textbf{\eref{condfeps}} with $\rho\geq 1$ and
then $\widehat{\theta}_1$ and $\widehat{\theta}_2$ are
consistent with a rate that depends on $\beta^0$.
Let us be more precise. According to the proof of Theorem \ref{thv1}, for
$j=1,\ldots,m+p$, the term $B_{n,j}^2(\theta^0)$ are of order $\exp(-2C_n/\beta^0)$ and the
term $V_{n,j}(\theta^0)/n$ are of order
$C_n^{2\alpha+(1-\rho)+(1-\rho_-)}\exp(-2C_n/\beta^0+2\delta C_n^\rho)/n.$

Set $C_n^*$ that realizes the best compromise between the squared bias
and the variance terms. It is independent from $\beta^0$ and  is given by 
$$C_n^*=\left[\frac{\log
n}{2\delta}-\frac{(2\alpha+(1-\rho)_-)}{2\delta\rho}\log\left(\frac{\log n}{2\delta}\right)\right]^{1/\rho}.$$
This choice yields to the rate
$$\varphi_n^2=\max\left\lbrace n^{-1}, \exp\left[-\frac{2}{\beta^0}\left(\frac{\log
n}{2\delta}-\frac{2\alpha+(1-\rho)_-}{2\delta \rho}\log\left(\frac{\log
n}{2\delta}\right)\right)^{1/\rho}\right](\log
n)^{(1-\rho)/\rho}\right\rbrace.$$ In other words, if $\rho=1$, then
$\mathbb{E}\parallel\widehat{\theta}_{1} -\theta^0\parallel_{\ell^2}^2
=O(1)\max\left\lbrace n^{-1}, n^{-1/(\beta^0 \delta)}(\log
n)^{2\alpha/(\beta^0 \delta)}\right\rbrace $  and if $\rho>1$, then
$\mathbb{E}\parallel\widehat{\theta}_{1} -\theta^0\parallel_{\ell^2}^2
=O(1)\exp\left[-{2}(\beta^0)^{-1}\left(\log
n/({2\delta})\right)^{1/\rho}\right].$

}

\end{exem}

\section{Comment on the use of the least squares criterion}
\setcounter{equation}{0}
\setcounter{lem}{0}
\setcounter{theo}{0}

In a proportional hazard model without errors, the main
drawback of the least squares criterion, compared to the partial log-likelihood, is that it does not allow to separate
the estimation of $\beta^0$ from the estimation of the baseline hazard
function $\eta$. The subject of this part is to motivate the choice of the
least squares criterion when the covariate is mismeasured.

First, as it is mentionned in the introduction, the partial
log-likelihood related to the filtration given by the observations only has not
an explicit form.

Second, consider as a partial log-likelihood related to the failure hazard 
function defined by \eref{induite} 
\begin{eqnarray*}
\frac{1}{n}\sum_{i=1}^n\int_0^\tau
\log\left\{ \frac{\mathbb{E}(f_{\beta}(Z_i)|\sigma(U_i,\ind_{T_i\geq t})}{n^{-1}\sum_{j=1}^n Y_j(t)
\mathbb{E}[f_{\beta}(Z_j)|\sigma(U_j,\ind_{T_j\geq t})]}\right\}dN_i(t).
\end{eqnarray*}
This partial log-likelihood depends on the observations, on the density $g$
of $Z$ and on $\eta_{\gamma^0}$, through $\mathbb{E}[f_{\beta}(Z)|\sigma(U,\ind_{T\geq t})]$. 
Hence the estimation of $\beta^0$
also depends on $\eta$ 
through the conditionning. 

Lastly, since the $Z_i$'s are unobservable, one other idea would be to
estimate $\beta^0$ by minimizing $\hat
L_n^{(1)}(\beta,U^{(n)})$ given by
\begin{multline*}
\hat
L_n^{(1)}(\beta,U^{(n)})=\frac{1}{n}\sum_{i=1}^n\int_0^\tau\left[
  \left(\frac{f_{\beta}^{(1)}W}{f_{\beta}}\right)\star K_{n,C_n}(U_i)\right.\\\left.-\frac{
    \sum_{j=1}^nY_j(t)(f_{\beta}^{(1)}W)\star K_{n,C_n}(U_j)}
{\sum_{j=1}^n
 Y_j(t)(f_{\beta} W) \star K_{n,C_n}(U_j)}\right]W\star K_{n,C_n}(U_i)dN_i(t), 
\end{multline*}
for $W$ a suitable chosen weight function, with $W(z)\not=0$ for all $z$ in
$\mathbb{R}$. Due to the
unobservability of $Z^{(n)}$,  $\hat
L_n^{(1)}(\beta,U^{(n)})$ can be seen as an estimation of the expectation of \eref{vraispart}.
Under reasonnable assumptions, $\hat\beta^{PL}$ such that
$\hat
L_n^{(1)}(\hat\beta^{PL},U^{(n)})=0$ is a consistent estimator of $\beta^0$. The main
difficulty lies in the study of its rate of
convergence. As in the study of
$\widehat{\theta}_1$, the rate of convergence of $\hat\beta^{PL}$ depends on
the smoothness of
${(f_{\beta}^{(1)}W)(z)}/{(f_{\beta})(z)}$, as a function of $z$, through the behavior of the ratio
$$\frac{\left({(f_{\beta^0}^{(1)}W})/{f_{\beta^0}}\right)^*(t)}{\overline{f_\varepsilon^*}(t)
}, \mbox{ as }t \rightarrow \infty.$$
Consequently, the best properties would be obtained for $W$ such that 
$(f_{\beta}^{(1)}W)/{(f_{\beta})}$ is in $\mathbb{L}_1(\mathbb{R})$ and has the
best smoothness properties. 
In the Cox model, $f^{(1)}_\beta(z)/f_\beta(z)=z$ and this estimation criterion 
provides $\sqrt{n}$-consistency and asymptotically Gaussian estimator, analogously
to the Nakamura's \citeyear{NAKAMURA92}'s estimator. The same result holds for
the relative risks considered  in Examples \ref{lemcos} and
\ref{lemcos2}. 
Nevertheless, for general relative risks, this criterion is less tractable than
the least squares criterion \eref{Sn1}, since it is strictly more difficult
to "smooth" $z\mapsto f^{(1)}_\beta(z)/f_\beta(z)$ than $z\mapsto f_\beta(z)$. This appears in a crucial
way in the model of excess relative risk where $f^{(1)}_\beta(z)/f_\beta(z)=z/(1+\beta
z)$.  This point has to be related to the difficulty and even the impossibility to find
a suitable correction of $L_n^{(1)}(\beta, U^{(n)})$, which leads to asymptotically
unbiased score functions (see the introduction).

\section{Proofs}
\setcounter{equation}{0}
\setcounter{lem}{0}
\setcounter{theo}{0}
From now $C$ denotes any numerical constant and  $C(A)$ indicates that it
depends on a $A$.
\subsection{Proof of Theorem \ref{thv1} }
\subsubsection{\textbf{Consistency}}
By classical arguments, the consistency  follows
from the two points :

\textbf{1-} The quantity $S_{\theta^0,g}(\theta)$ is minimum if and only if $\theta=\theta^0.$

\textbf{2-}
For all $\theta \in \Theta$, $\mathbb {E}[S_{n,1}(\theta) -
S_{\theta^0,g}(\theta)]^2=o(1)$ as $\nti $, with
$S_{\theta^0,g}(\theta)$ defined in \eref{S0},

\textbf{3-} If  ${\omega}(n,\rho)$ denotes 
${\omega}(n,\rho)=\sup\left\lbrace \vert {S_{n,1}} (\theta)-S_{n,1}(\theta^\prime)\vert:\|
\theta-\theta^\prime\|_{\ell^2}\leq \rho \right\rbrace,$ there exists $\rho_k$  tending to 0, such
that $\mathbb{E}[{\omega}(n,\rho_k)]^2=O(\rho_k^2)\mbox{ as }\nti ~~~~~~\forall k \in\mathbb{N}.$

\textbf{Proof of 1-}
Under \eref{maxunique}, by applying \eref{carre} we get 
\begin{eqnarray*}
\frac{\partial}{\partial
  \beta}S_{\theta^0,g}(\theta)=2\int_0^\tau 
\mathbb{E}\left[f_{\beta}^{(1)}(Z_i)\eta_{\gamma}(t)\left\lbrace\eta_{\gamma}(t)f_{\beta}(Z)-\eta_{\gamma^0}(t)f_{\beta^0}(Z)\right\rbrace
W(Z)Y(t)\right]dt=0
\Leftrightarrow \theta=\theta^0,
\end{eqnarray*}
and
\begin{eqnarray*}
\frac{\partial}{\partial
  \gamma}S_{\theta^0,g}(\theta)=2\int_0^\tau 
\mathbb{E}\left[f_{\beta}(Z_i)\eta^{(1)}_{\gamma}(t)\left\lbrace\eta_{\gamma}(t)f_{\beta}(Z)-\eta_{\gamma^0}(t)f_{\beta^0}(Z)\right\rbrace
W(Z)Y(t)\right]dt=0
\Leftrightarrow \theta=\theta^0.
\end{eqnarray*}
The matrix of second derivatives equals
\begin{eqnarray*}
\left(\frac{\partial^2 S_{\theta^0,g}(\theta)}{\partial
  \theta^2}\left.\right\vert_{\theta=\theta^0}\right)
&:=&H(\theta^0)=\begin{pmatrix} H_{11}(\theta^0)& H_{12}(\theta^0)\\
(H_{12}(\theta^0))^\top & H_{22}(\theta^0)
\end{pmatrix}
\end{eqnarray*}
with
\begin{eqnarray*}
H_{11}(\theta^0)
&=&2\int_0^\tau
  \mathbb{E}\left[(f_{\beta^0}^{(1)}(Z))(f_{\beta^0}^{(1)}(Z))^\top\eta_{\gamma^0}^2(t)W(Z)Y(t)\right]dt\\
H_{12}(\theta^0)
&=&2  \int_0^\tau
\mathbb{E}\left[f_{\beta^0}(Z)\eta_{\gamma^0}(t)f_{\beta^0}^{(1)}(Z)(\eta_{\gamma^0}^{(1)}(t))^\top
W(Z)Y(t)\right]dt\\
H_{22}(\theta^0)
&=& 2 \int_0^\tau
  \mathbb{E}\left[(\eta_{\gamma^0}^{(1)}(t))(\eta_{\gamma^0}^{(1)}(t))^\top f_{\beta^0}^2(Z)W(Z)Y(t)\right]dt.
\end{eqnarray*}
An obvious application of Cauchy-Schwarz Inequality gives that the matrix $H$
is non negative definite and hence under \eref{concls}
\textbf{1-} is proved.

\textbf{Proof of 2-}

For both the bias and the variance, we will give two upper bounds, based on
the two following applications of the
H\"older's inequality
\begin{eqnarray}
\label{holder1}
\vert<\varphi_1,\varphi_2>\vert\leq \parallel
\varphi_1\parallel_2\parallel\varphi_2\parallel_2,
\end{eqnarray}
 and 
\begin{eqnarray}\label{holder2}
\vert <\varphi_1,\varphi_2>\vert\leq 
\parallel \varphi_1\parallel_\infty\parallel\varphi_2\parallel_1.\end{eqnarray}
According to Lemma \ref{lemdeconv} we write that
\begin{eqnarray*}
\mathbb{E}[S_{n,1}(\theta)]
=\int_0^\tau \mathbb{E}\left[(f_{\beta}^2W)\star
  K_{C_n}(Z)\eta^2_{\gamma}(t)Y(t)-2(f_{\beta} W)\star
  K_{C_n}(Z)\eta_{\gamma}(t)f_{\beta^0}(Z)\eta_{\gamma^0}(t)Y(t)\right] dt,
\end{eqnarray*}
and hence
\begin{eqnarray*}
\mathbb{E}(S_{n,1}(\theta))-S_{\theta^0,g}(\theta)
&=&\int\int_0^\tau \eta_{\gamma}^2(t)\ind_{x\geq t}\left<f_{X,Z}(x,\cdot),(f_{\beta}^2W)\star
  K_{C_n}-f_{\beta}^2W\right>dx\,dt\\&&-2\int\int_0^\tau
\eta_{\gamma^0}(t)\eta_{\gamma}(t)\ind_{x\geq t}\left<f_{\beta^0}(\cdot)f_{X,Z}(x,\cdot),(f_{\beta} W)\star
  K_{C_n}-f_{\beta}
  W \right>dx\,dt.
\end{eqnarray*}
By applying  \eref{holder1} we obtain the first bound
\begin{eqnarray*}
\left\vert\mathbb{E}(S_{n,1}(\theta))-S_{\theta^0,g}(\theta)\right\vert&\leq&
\Big(\int_0^\tau\eta_\gamma^2(t)dt\Big) \parallel f_{X,Z}\parallel_2\parallel (f_\beta^2
W)\star K_{C_n}-(f_\beta^2 W)\parallel_2\\&&+
\Big(2\int_0^\tau\eta_\gamma(t)\eta_{\gamma^0}(t)dt\Big) \int \parallel
f_{\beta^0}(\cdot)f_{X,Z}(x,\cdot)\parallel_2 dx\parallel (f_\beta
W)\star K_{C_n}-f_\beta W\parallel_2.
\end{eqnarray*}
Applying Parseval's formula we get
\begin{eqnarray*}
\left\vert\mathbb{E}(S_{n,1}(\theta))-S_{\theta^0,g}(\theta)\right\vert
&\leq& (2\pi)^{-1}\Big(\int_0^\tau\eta_\gamma^2(t)dt\Big) \parallel f_{X,Z}\parallel_2\parallel (f_\beta^2
W)^*(K_{C_n}^*-1)\parallel_2
\\&&+(\pi)^{-1}\Big(\int_0^\tau\eta_\gamma(t)\eta_{\gamma^0}(t)dt\Big) \int \parallel
f_{\beta^0}(\cdot)f_{X,Z}(x,\cdot)\parallel_2 dx
\parallel (f_\beta
W)^*(K_{C_n}^*-1)\parallel_2
\end{eqnarray*}
that is
\begin{multline}
\label{bc1}
\left\vert\mathbb{E}(S_{n,1}(\theta))-S_{\theta^0,g}(\theta)\right\vert\leq C(\gamma,\gamma^0,f_{\beta^0})
\\\left[\parallel(f_\beta^2
W)^*(K_{C_n}^*-1)\parallel_2+\parallel(f_\beta W)^*(K_{C_n}^*-1)\parallel_2\right].
\end{multline}
According to \eref{holder2},  $\left\vert
\mathbb{E}(S_{n,1}(\theta))-S_{\theta^0,g}(\theta)\right\vert$ is also bounded by
\begin{multline*}
\parallel
(f_{\beta}^2W)\star
K_{C_n}-(f_{\beta}^2W)\parallel_\infty  \int \parallel
f_{X,Z}(x,\cdot)\parallel_1 dx\Big(\int_0^\tau \eta^2_{\gamma}(t)
dt\Big)\\+
\parallel (f_{\beta} W)\star K_{C_n}-(f_{\beta} W)\parallel_\infty \Big(
2 \int \parallel f_{\beta^0}(\cdot)f_{X,Z}(x,\cdot)\parallel_1dx\int_0^\tau\eta_{\gamma}(t)\eta_{\gamma^0}(t)
dt\Big)\\
\leq \parallel
(f_{\beta}^2W)^*
(K_{C_n}^*-1)\parallel_1\Big((2\pi)^{-1}\parallel f_{X,Z}\parallel_1\int_0^\tau\eta^2_{\gamma}(t)
dt\Big)\\+
\parallel (f_{\beta} W)^*(K_{C_n}^*-1)\parallel_1
\Big(\pi^{-1} \int \parallel f_{\beta^0}(\cdot)f_{X,Z}(x,\cdot)\parallel_1 dx\int_0^\tau\eta_{\gamma}(t)f_{\beta^0}\eta_{\gamma^0}(t)
dt\Big).
\end{multline*}
This implies that
\begin{eqnarray*}
\left\vert\mathbb{E}(S_{n,1}(\theta))-S_{\theta^0,g}(\theta)\right \vert&\leq&
\Big[\int_0^\tau
\eta_{\gamma}^2(t)dt\Big]\parallel(f_{\beta}^2W)^*(K_{C_n}^*-1)\parallel_1\\&&+\Big[\mathbb{E}\vert
  f_{\beta^0}(Z)\vert\int_0^\tau
\eta_{\gamma}(t)\eta_{\gamma^0}(t)dt\Big]\parallel(f_{\beta} W)^*(K_{C_n}^*-1)\parallel_1,
\end{eqnarray*}
that is
\begin{multline}
\label{bc2}
\left\vert\mathbb{E}(S_{n,1}(\theta))-S_{\theta^0,g}(\theta)\right \vert\leq
C(\gamma,\gamma^0,f_{\beta^0})\\\left[\parallel(f_{\beta} W)^*(K_{C_n}^*-1)\parallel_1
 +\parallel(f_\beta^2W)^*(K_{C_n}^*-1)\parallel_1\right].
\end{multline}
By combining the bounds \eref{bc1} and \eref{bc2}  we get that
\begin{multline}
\label{bc}
\left\vert\mathbb{E}(S_{n,1}(\theta))-S_{\theta^0,g}(\theta)\right \vert\leq
C(\gamma,\gamma^0,f_{\beta^0})\\\times\min\left\lbrace
\parallel(f_{\beta}
W)^*(K_{C_n}^*-1)\parallel_2+\parallel(f_\beta^2W)^*(K_{C_n}^*-1)\parallel_2,\right.\\\left.\parallel(f_{\beta}
W)^*(K_{C_n}^*-1)\parallel_1+\parallel(f_\beta^2W)^*(K_{C_n}^*-1)\parallel_1\right\rbrace.\end{multline}
By applying Lemma \ref{contint}
$$
\left\vert\mathbb{E}(S_{n,1}(\theta))-S_{\theta^0,g}(\theta)\right \vert^2
=O\Big(C_n^{-2a+1-r+(1-r)_-}\exp(-2dC_n^r)\Big)=o(1).
$$
\textbf{Study of the variance}
Since the random variables are i.i.d., we get that
\begin{eqnarray*}
\mbox{Var}[S_{n,1}(\theta)]
&=&\frac{(2+o(1))}{n}(A_1+A_2),
\end{eqnarray*}
with
\begin{eqnarray*}
A_1&=&\mathbb{E}\Big[
  (f_{\beta}^2 W)\star K_{n,C_n}(U)\int_0^\tau\eta_{\gamma}^2(t)Y(t)dt\Big]^2\\
\mbox{ and }A_2&=&4\mathbb{E}\Big[
  (f_{\beta} W)\star K_{n,C_n}(U)\int_0^\tau\eta_{\gamma}(t)dN(t)\Big]^2.
\end{eqnarray*}
According to \eref{holder2} and  by applying Lemma \ref{lemdeconv}, $A_1$ is less than
\begin{multline*}
\Big(\int_0^\tau\eta_{\gamma}^2(t)dt\Big)^2 \int
\left\vert\left< f_{X,Z}(x,\cdot)\star f_\varepsilon,((f_\beta^2 W)\star
K_{n,C_n})^2\right>\right\vert dx\\\leq 
\Big(\int_0^\tau\eta_{\gamma}^2(t)dt\Big)^2 \int \parallel
f_{X,Z}(x,\cdot)\star f_\varepsilon\parallel_\infty dx
\parallel (f_\beta^2 W)\star K_{n,C_n}\parallel_2^2
\end{multline*}
and hence
\begin{eqnarray*}
A_1
&\leq& (2\pi)^{-1}\Big(\int_0^\tau\eta_{\gamma}^2(t)dt\Big)^2\parallel
f_\varepsilon\parallel_\infty
\parallel f_{X,Z}\parallel_1\left\|\frac{(f_{\beta}^2W)^*K_{C_n}^*}{\overline{f_\varepsilon^*}} \right\|_2^2.
\end{eqnarray*}
In the same way, we get a first bound for $A_2$.
Let us denote by
\begin{eqnarray}
\label{defphi}
\varphi(X,Z)=\int_0^\tau \eta_\gamma(t)dN(t).
\end{eqnarray}
According to Lemma \ref{lemdeconv} and to \eref{holder2}, $A_2$ is bounded by
\begin{multline*}
4 \int \left<(\varphi^2(x,\cdot)f_{X,Z}(x,\cdot))\star f_\varepsilon, ((f_\beta
W)\star K_{n,C_n})^2 \right>
dx\\
\leq 4\int \parallel (\varphi^2(x,\cdot)f_{X,Z}(x,\cdot))\star
f_\varepsilon\parallel_\infty dx \parallel(f_\beta
W)\star K_{n,C_n}\parallel_2^2.
\end{multline*}
Since 
$$ \int \parallel \varphi^2(x,\cdot)f_{X,Z}(x,\cdot)\parallel_1dx=\mathbb{E}\Big[\int_0^\tau \eta_\gamma(t)dN(t)\Big]^2,$$
we get that
$$ \int \parallel \varphi^2(x,\cdot)f_{X,Z}(x,\cdot))\star
f_\varepsilon\parallel_\infty dx\leq \parallel f_\varepsilon\parallel_\infty\mathbb{E}\Big[\int_0^\tau \eta_\gamma(t)dN(t)\Big]^2
.$$
Consequently,
\begin{eqnarray*}
A_2&\leq& 4\left[(2\pi)^{-1}\mathbb{E}\Big(\int_0^\tau\eta_{\gamma}(t)dN(t)\Big)^2\parallel
f_\varepsilon\parallel_\infty \right]\left\|\frac{(f_{\beta}
    W)^*K_{C_n}^*}{\overline{f_\varepsilon^*}} \right\|_2^2.
\end{eqnarray*}
It follows that, 
\begin{eqnarray}
\label{bva1}
\qquad\mbox{Var}[S_{n,1}(\theta)]\leq
\frac{C(\theta^0,\parallel f_\varepsilon\parallel_\infty)}{n}\Big[\Big\|(f_{\beta}
    W)^*\frac{K_{C_n}^*}{\overline{f_\varepsilon^*}} \Big\|_2^2+\Big\|
(f_{\beta}^2W)^*\frac{K_{C_n}^*}{\overline{f_\varepsilon^*}} \Big\|_2^2\Big].
\end{eqnarray}
According to \eref{holder2},  $A_1$ is also less
than
% echange des roles entre phi_1 et phi_2
\begin{multline*}
\Big(\int_0^\tau\!\!\!\eta_{\gamma}^2(t)dt\Big)^2\int \left\vert
\left< f_{X,Z}(x,\cdot)\star f_\varepsilon, ((f_\beta^2W)\star
K_{n,C_n})^2\right>\right\vert dx\\\leq 
\Big(\int_0^\tau\!\!\!\eta_{\gamma}^2(t)dt\Big)^2\parallel \int
f_{X,Z}(x,\cdot)\parallel_1 dx\parallel (f_\beta^2 W)\star K_{n,C_n}\parallel_\infty
^2.\end{multline*}
In the same way $A_2$ is less than
\begin{multline*}
 4\int \left<(\varphi^2(x,\cdot)f_{X,Z}(x,\cdot))\star f_\varepsilon, ((f_\beta
W)\star K_{n,C_n})^2 \right>
dx\\
\leq 4\int \parallel (\varphi^2(x,\cdot)f_{X,Z}(x,\cdot))\star
f_\varepsilon\parallel_1dx \parallel(f_\beta
W)\star K_{n,C_n}\parallel_\infty^2
\end{multline*}
where $\varphi(X,Z)$ is defined in \eref{defphi}.
Once again, since
$$ \int \parallel (\varphi^2(x,\cdot)f_{X,Z}(x,\cdot))\star
f_\varepsilon\parallel_1dx=\mathbb{E}\Big(\int_0^\tau \eta_\gamma(t)dN(t)\Big)^2
,$$
\begin{eqnarray}
\label{bva2}
\mbox{Var}[S_{n,1}(\theta)]\leq \frac{C(\theta^0) }{n} 
\Big[\Big\|(f_\beta
W)^*\frac{K_{C_n}^*}{\overline{f_\varepsilon^*}}\Big\|_1^2+\Big\|(f_\beta^2 W)^* 
\frac{K_{C_n}^*}{\overline{f_\varepsilon^*}}\Big\|_1^2\Big].
\end{eqnarray}
By combining \eref{bva1} and \eref{bva2}, we obtain that
\begin{multline}
\label{bva}
\mbox{Var}[S_{n,1}(\theta)] \leq \frac{C(\theta^0, \parallel f_\varepsilon\parallel_\infty)}{n} 
\min\left\lbrace \Big\|(f_\beta
W)^*\frac{K_{C_n}^*}{\overline{f_\varepsilon^*}}\Big\|_2^2+\Big\|(f_\beta^2 W)^*
\frac{K_{C_n}^*}{\overline{f_\varepsilon^*}}\Big\|_2^2,\right.\\ \left.\qquad\Big\|(f_\beta
W)^*\frac{K_{C_n}^*}{\overline{f_\varepsilon^*}}\Big\|_1^2+\Big\|(f_\beta^2 W)^*
\frac{K_{C_n}^*}{\overline{f_\varepsilon^*}}\Big\|_1^2\right\rbrace.\end{multline}
According to Lemma \ref{contint}
$$\mbox{Var}[S_{n,1}(\theta)]
=O\Big(C_n^{2(\alpha-a)+1-\rho+(1-\rho)_-}\exp(-2dC_n^r+2\delta
C_n^\rho)/n\Big),$$
and hence under \eref{condcons2},
$\mathbb{E}\left[S_{n,1}(\theta)-S_{\theta^0,g}(\theta)\right]^2=o(1),$ as $n\rightarrow \infty.$

\textbf{Proof of 3-}

By definition $S_{n,1}(\theta)-S_{n,1}(\theta^\prime)$ equals
\begin{multline*}
-\frac{2}{n}\int_0^\tau\left[(f_{\beta} W)\star
K_{n,C_n}(U_i)\eta_{\gamma}(t)-(f_{\beta^\prime} W)\star
K_{n,C_n}(U_i)\eta_{\gamma^\prime}(t)\right]dN_i(t)\\+
\frac{1}{n}\int_0^\tau\left[(f_{\beta}^2 W)\star
K_{n,C_n}(U_i)\eta^2_{\gamma}(t)-(f_{\beta^\prime}^2 W)\star
K_{n,C_n}(U_i)\eta^2_{\gamma^\prime}(t)\right]Y_i(t)dt.
\end{multline*}
Under \textbf{\eref{deriv}}, \textbf{\eref{cw1}}, \textbf{\eref{cw11}}, \textbf{\eref{condfeps}} and \textbf{\eref{super}}, for $C_n$ satisfying \eref{condcons2}, since $\|\theta-\theta^\prime\|_{\ell^2}\leq \rho_k$,
we get that $\mathbb{E}(\vert
S_{n,1}(\theta)-S_{n,1}(\theta^\prime)\vert^2)=O(\rho_k^2)$. Hence
\textbf{3-} follows. \hfill $\Box$

\subsubsection{\bf Rate of convergence}
Denote by $S_{n,1}^{(1)}(\theta)$ and $S_{n,1}^{(2)}(\theta)$ the first and
second derivatives of $S_{n,1}(\theta)$ with respect to $\theta$. By using
classical Taylor 
expansion 
and the consistency of $\widehat{\theta}_1$, we get that
$
0=S^{(1)}_{n,1}(\widehat{\theta}_1)=S^{(1)}_{n,1}(\theta^0)+
S^{(2)}_{n,1}(\theta^0)(\widehat{\theta}_1-\theta^0)+
R_{n}(\widehat{\theta}_1-\theta^0),$
 with $R_{n}$ defined by
\begin{eqnarray}
R_{n}=\int_{0}^{1}[S^{(2)}_{n,1}(\theta^0+s(\widehat{\theta}_1-\theta^0))
-S^{(2)}_{n,1}(\theta^0)]ds. \label{Rn}
\end{eqnarray}
This implies that
\begin{eqnarray}
\label{base} \widehat{\theta}_1-\theta^0=-[S^{(2)}_{n,1}(\theta^0)+R_{n}]^{-1}
S^{(1)}_{n,1}(\theta^0).\end{eqnarray} 
Consequently we have to check the three following points
\begin{itemize}
\item[i)]
  $\mathbb{E}\Big[\big\{S^{(1)}_{n,1}(\theta^0))-S^{(1)}_{\theta^0,g}(\theta^0)\big\}\big\{S^{(1)}_{n,1}(\theta^0))-
S^{(1)}_{\theta^0,g}(\theta^0)\big\}^\top\Big]=
O[\varphi_n\varphi_n^\top]$, 
\item[ii)]
$\mathbb{E}\left[S^{(2)}_{n,1}(\theta^0)-
  S^{(2)}_{\theta^0,g}(\theta^0)\right]^2=o(1)$, 
\item[iii)] $R_{n}$ defined in \eref{Rn} satisfies
  $\mathbb{E}(\parallel R_{n}\parallel_{\ell^2}^2)=o(1)$ as $n\rightarrow \infty$.
\item[iv)]
$\mathbb{E}\|\widehat{\theta}_1-\theta^0\|_{\ell^2}^2\leq
4\mathbb{E}\left[
(S^{(1)}_{n,1}(\theta^0))^\top\Big[\Big(\frac{\partial^2 S_{\theta^0,g}(\theta)}{\partial
\theta_j\theta_k}|_{\theta=\theta^0}\Big)^{-1}\Big]
^{\top}\Big(\frac{\partial^2 S_{\theta^0,g}(\theta)}{\partial
\theta_j\theta_k}|_{\theta=\theta^0}\Big)^{-1} S^{(1)}_{n,1}(\theta^0)\right]+o(\varphi_{n}^2).$
\end{itemize}
The rate of convergence of $\widehat{\theta}_1$ is thus given by the
order of $S^{(1)}_{n,1}(\theta^0)-S^{(1)}_{\theta^0,g}(\theta^0)=S^{(1)}_{n,1}(\theta^0)$. 
\begin{center}
\textbf{Proof of i)}
\end{center}
According to \eref{Sn1}, $S^{(1)}_{n,1}(\theta^0)$ equals
\begin{multline}
\label{Sn1b}
\frac{2}{n}\sum_{i=1}^n\begin{pmatrix}\displaystyle
-\int_0^\tau (f^{(1)}_{\beta^0} W)\star
    K_{n,C_n}(U_i)\eta_{\gamma^0}(t)dN_i(t)+\int_0^\tau (f_{\beta^0} f_{\beta^0}^{(1)}W)\star
    K_{n,C_n}(U_i)\eta_{\gamma^0}^2(t)Y_i(t)dt \\\\\displaystyle
-\int_0^\tau (f_{\beta^0} W)\star
    K_{n,C_n}(U_i)\eta^{(1)}_{\gamma^0}(t)dN_i(t)+\int_0^\tau (f_{\beta^0}^2W)\star
    K_{n,C_n}(U_i)\eta_{\gamma^0}(t)\eta_{\gamma^0}^{(1)}(t)Y_i(t)dt
\end{pmatrix}.
\end{multline}
\textbf{Study of the bias}
By definition, $\mathbb{E}({\partial
S_{n,1}(\theta)}/{\partial \beta})_{\theta=\theta^0}$ equals
$$ -2\mathbb{E}\Big[\int_0^\tau (f^{(1)}_{\beta^0} W)\star
    K_{n,C_n}(U_1)\eta_{\gamma^0}(t)dN_1(t)\Big]+2\mathbb{E}\Big[\int_0^\tau (f_{\beta^0} f_{\beta^0}^{(1)}W)\star
    K_{n,C_n}(U_1)\eta_{\gamma^0}^2(t)Y_1(t)dt\Big].$$
Hence, according to Lemma \ref{lemdeconv},
\begin{eqnarray*}
\mathbb{E}\left(\frac{\partial
S_{n,1}(\theta)}{\partial \beta}\left.\right\vert_{\theta=\theta^0}\right)&=&-2\mathbb{E}\left[f_{\beta^0}(Z_1) (f^{(1)}_{\beta^0} W)\star
    K_{n,C_n}(U_1)\int_0^\tau
\eta_{\gamma^0}^2(t)Y_1(t)dt\right]\\
&&+2\mathbb{E}\left[(f_{\beta^0} f_{\beta^0}^{(1)}W)\star
    K_{n,C_n}(U_1)\int_0^\tau\eta_{\gamma^0}^2(t)Y_1(t)dt\right]\\
&=&-2\mathbb{E}\left[f_{\beta^0}(Z_1) (f^{(1)}_{\beta^0} W)\star
    K_{C_n}(Z_1)\int_0^\tau
\eta_{\gamma^0}^2(t)Y_1(t)dt\right]
\\&&+2\mathbb{E}\left[(f_{\beta^0} f_{\beta^0}^{(1)}W)\star
    K_{C_n}(Z_1)\int_0^\tau\eta_{\gamma^0}^2(t)Y_1(t)dt\right].
\end{eqnarray*}
Since 
$$\left.\frac{\partial S_{\theta^0,g}(\theta)}{\partial
\beta}\right\vert_{\theta=\theta^0}\!\!\!\!\!\!\!=
-2\mathbb{E}\left[
f_{\beta^0}(Z_1)(f^{(1)}_{\beta^0} W)(Z_1)\int_0^\tau \eta_{\gamma^0}^2(t)Y_1(t)dt\right]
+2\mathbb{E}\left[(f_{\beta^0}
f_{\beta^0}^{(1)}W)(Z_1)\int_0^\tau\eta_{\gamma^0}^2(t)Y_1(t)dt\right]=0,$$ we get that 
 $\mathbb{E}\left(\left.\partial S_{n,1}(\theta)/\partial
\beta\right\vert_{\theta=\theta^0}\right)$ also equals
\begin{multline*}
-2\mathbb{E}\left[f_{\beta^0}(Z_1) [(f^{(1)}_{\beta^0} W)\star
    K_{C_n}(Z_1)-(f^{(1)}_{\beta^0} W)(Z_1)]\int_0^\tau \eta_{\gamma^0}^2(t)Y_1(t)dt\right]\\+2\mathbb{E}\left[[(f_{\beta^0} f_{\beta^0}^{(1)}W)\star
    K_{C_n}(Z_1)-(f_{\beta^0}
f_{\beta^0}^{(1)}W)(Z_1)]\int_0^\tau\eta_{\gamma^0}^2(t)Y_1(t)dt\right]\\
=-2\int \left<f_{\beta^0}(\cdot)f_{X,Z}(x,\cdot), [(f^{(1)}_{\beta^0} W)\star
    K_{C_n}-(f^{(1)}_{\beta^0} W)]\right>\int_0^\tau
\eta_{\gamma^0}^2(t)\ind _{x\geq t}dt\,dx\\+2\int \left<f_{X,Z}(x,\cdot),[(f_{\beta^0} f_{\beta^0}^{(1)}W)\star
    K_{C_n}-(f_{\beta^0}
f_{\beta^0}^{(1)}W)]\right>\int_0^\tau\eta_{\gamma^0}^2(t)\ind_{x\geq
t}dt\,dx.
\end{multline*}
In the same way, according to Lemma \ref{lemdeconv},
\begin{eqnarray*}
\mathbb{E}\left(\frac{\partial S_{n,1}(\theta)}{\partial \gamma}\left.\right\vert_{\theta=\theta^0}\right)&=&2\mathbb{E}\left[- (f_{\beta^0} W)\star
    K_{n,C_n}(U_1)\int_0^\tau\eta^{(1)}_{\gamma^0}(t)dN_1(t)\right]\\&&+ 2\mathbb{E}\left[(f_{\beta^0}^2W)\star
    K_{n,C_n}(U_1)\int_0^\tau\eta_{\gamma^0}(t)\eta_{\gamma^0}(t)Y_1(t)dt\right]\\
&=&-2\mathbb{E}\left[f_{\beta^0}(Z_1)(f_{\beta^0}W)\star
    K_{n,C_n}(U_1)\int_0^\tau \eta^{(1)}_{\gamma^0}(t)
\eta_{\gamma^0}(t)Y_1(t)dt\right]\\&&+2\mathbb{E}\left[(f_{\beta^0}^2W)\star
    K_{n,C_n}(U_1)\int_0^\tau \eta_{\gamma^0}^{(1)}(t)\eta_{\gamma^0}(t)Y_1(t)dt\right]\\
&=&-2\mathbb{E}\left[f_{\beta^0}(Z_1) (f_{\beta^0} W)\star
    K_{C_n}(Z_1)\int_0^\tau \eta^{(1)}_{\gamma^0}(t)\eta_{\gamma^0}(t)Y_1(t)dt\right]\\&&+2\mathbb{E}\left[(f_{\beta^0}^2W)\star
    K_{C_n}(Z_1)\int_0^\tau \eta_{\gamma^0}^{(1)}(t)\eta_{\gamma^0}(t)Y_1(t)dt\right].
\end{eqnarray*}
Since 
\begin{eqnarray*}
\frac{\partial S_{\theta^0,g}(\theta)}{\partial
\gamma_j}\left.\right\vert_{\theta=\theta^0}&=&
-2\mathbb{E}\left[f_{\beta^0}(Z_1)(f_{\beta^0} W)(Z_1)\int_0^\tau
\eta_{\gamma^0}^{(1)}\eta_{\gamma^0}(t)Y_1(t)dt\right]\\&&+2\mathbb{E}\left[(f_{\beta^0}^2W)(Z_1)\int_0^\tau\eta_{\gamma^0}^{(1)}(t)\eta_{\gamma^0}(t)Y_1(t)dt\right]\\&=&0,\end{eqnarray*}
 we obtain that $\mathbb{E}\left({\partial S_{n,1}(\theta)}/{\partial
\gamma}\left.\right\vert_{\theta=\theta^0}\right)$ equals
\begin{eqnarray*}
&-&2\mathbb{E}\left[f_{\beta^0}(Z_1) [(f_{\beta^0} W)\star
    K_{C_n}(Z_1)-(f_{\beta^0} W)(Z_1)]\int_0^\tau \eta_{\gamma^0}^{(1)}\eta_{\gamma^0}(t)Y_1(t)dt\right]\\&&+2\mathbb{E}\left[[(f_{\beta^0}^2W)\star
    K_{C_n}(Z_1)-(f_{\beta^0}^2W)(Z_1)]\int_0^\tau\eta_{\gamma^0}^{(1)}(t)\eta_{\gamma^0}(t)Y_1(t)dt\right]\\
\!\!\!&=&\!\!\!-2\int \left<f_{\beta^0}(\cdot)f_{X,Z}(x,\cdot), [(f_{\beta^0} W)\star
    K_{C_n}-(f_{\beta^0} W)]\right>\int_0^\tau 
\eta^{(1)}_{\gamma^0}(t)\eta_{\gamma^0}(t)\ind _{x\geq t}dt\,dx\\&&+2\int \left<f_{X,Z}(x,\cdot),[(f_{\beta^0}^2W)\star
    K_{C_n}-(f_{\beta^0}^2W)]\right>\int_0^\tau \eta_{\gamma^0}^{(1)}(t)\eta_{\gamma^0}(t)\ind_{x\geq
t}dt\,dx.
\end{eqnarray*}
A first bound for this bias term can be obtained by writing that for $j=1,\cdots,m$
$$(1/2)\mathbb{E}\left(\left.\partial S_{n,1}(\theta)/\partial
\beta_j\right\vert_{\theta=\theta^0}\right)
$$ is bounded by
\begin{multline*}
\parallel
\left( f^{(1)}_{\beta^0,j} W\right)\star    K_{C_n}-\left(f^{(1)}_{\beta^0,j}
W\right)\parallel_2
\left(\int \parallel f_{\beta^0}(\cdot)f_{X,Z}(x,\cdot)\parallel_2dx 
\int_0^\tau
\eta_{\gamma^0}^2(t)dt\right)\\+ \left\| \left(f^{(1)}_{\beta^0,j}f_{\beta^0}W\right)\star
    K_{C_n}-\left(f^{(1)}_{\beta^0,j}f_{\beta^0}W\right)\right\|_2\left(\int \parallel
f_{X,Z}(x,\cdot)\parallel_2dx\int_0^\tau\eta_{\gamma^0}^2(t)dt\right)\\
\leq 
\left\| \left(f^{(1)}_{\beta^0,j} W\right)^*
  (K_{C_n}^*-1)\right\|_2\left((2\pi)^{-1}\int \parallel f_{\beta^0}(\cdot)f_{X,Z}(x,\cdot)\parallel_2dx\int_0^\tau
\eta_{\gamma^0}^2(t)dt\right)\\+ \left\|\left( f^{(1)}_{\beta^0,j}f_{\beta^0} W\right)^*
   (K_{C_n}^*-1)\right\|_2\left((2\pi)^{-1}\int \parallel
f_{X,Z}(x,\cdot)\parallel_2dx\int_0^\tau\eta_{\gamma^0}^2(t)dt\right).
 \end{multline*}
In the same way,
$$(1/2)\mathbb{E}\left(\left.\partial S_{n,1}(\theta)/\partial
\gamma_j\right\vert_{\theta=\theta^0}\right)
$$ is bounded by
\begin{multline*}
 \parallel
(f_{\beta^0} W)\star
    K_{C_n}-f_{\beta^0} W\parallel_2\left(\int \parallel f_{\beta^0}(\cdot)f_{X,Z}(x,\cdot)\parallel_2dx\int_0^\tau\left\vert
\eta^{(1)}_{\gamma^0,j}(t)\right\vert\eta_{\gamma^0}(t)dt\right)\\+ \parallel (f_{\beta^0}^2W)\star
    K_{C_n}-f_{\beta^0}^2 W\parallel_2\left(\int \parallel
f_{X,Z}(x,\cdot)\parallel_2dx\int_0^\tau\left\vert
\eta^{(1)}_{\gamma^0,j}(t)\eta_{\gamma^0}(t)\right\vert dt\right)\\
\leq 
\parallel (f_{\beta^0} W)^*
    (K_{C_n}^*-1)\parallel_2\left((2\pi)^{-1}\int \parallel f_{\beta^0}(\cdot)f_{X,Z}(x,\cdot)\parallel_2dx\int_0^\tau \left\vert 
\eta^{(1)}_{\gamma^0,j}(t)\right\vert 
\eta_{\gamma^0}(t)dt\right)\\+ \parallel (f_{\beta^0}^2 W)^*
    (K_{C_n}^*-1)\parallel_2\left((2\pi)^{-1}\int \parallel
f_{X,Z}(x,\cdot)\parallel_2dx\int_0^\tau\left\vert 
\eta^{(1)}_{\gamma^0,j}(t)\eta_{\gamma^0}(t)\right\vert dt\right).
 \end{multline*}
Consequently
\begin{multline}
\label{b1b1}
\left\vert \mathbb{E}\left(\frac{\partial S_{n,1}(\theta)}{\partial
\beta_j}\left.\right\vert_{\theta=\theta^0}\right)\right\vert
\\\leq C(\theta^0)\left[\left\|\left(f^{(1)}_{\beta^0,j}W\right)^*(K_{C_n}^*-1)\right\|_2+
\left\|\left(
(f^{(1)}_{\beta^0,j}f_{\beta^0}W\right)^*(K_{C_n}^*-1)\right\|_2\right],
\end{multline}
and
\begin{multline}
\label{b1b2}
\left\vert \mathbb{E}\left(\frac{\partial S_{n,1}(\theta)}{\partial
\gamma_j}\left.\right\vert_{\theta=\theta^0}\right)-\left(\frac{\partial S_{\theta^0,g}(\theta)}{\partial
\gamma_j}\left.\right\vert_{\theta=\theta^0}\right)\right\vert\\\leq
C(\theta^0)\left[\left\|\left(f_{\beta^0}W\right)^* (K_{C_n}^*-1)\right\|_2+\left\|\left(f_{\beta^0}^2 W
\right)^* (K_{C_n}^*-1)\right\|_2\right].
\end{multline}
A second bound for the bias term can be obtained by writing that
$(1/2)\mathbb{E}\left(\partial S_{n,1}(\theta)/\partial
\beta_j\left.\right\vert_{\theta=\theta^0}\right)$ is bounded by
\begin{multline*}
 \parallel
\left(f^{(1)}_{\beta^0,j} W\right)\star
    K_{C_n}-\left(f^{(1)}_{\beta^0,j} W\right)\parallel_\infty \left(\int \parallel f_{\beta^0}(\cdot)f_{X,Z}(x,\cdot)\parallel_1 dx\int_0^\tau
\eta_{\gamma^0}^2(t)dt\right)\\+ \left\| \left(f^{(1)}_{\beta^0,j}f_{\beta^0} W\right)\star
    K_{C_n}-\left(f^{(1)}_{\beta^0,j}f_{\beta^0} W\right)
\right\|_\infty\left(\int \parallel
f_{X,Z}(x,\cdot)\parallel_1 dx \int_0^\tau\eta_{\gamma^0}^2(t)dt\right)\\
\leq  \left\| \left(f^{(1)}_{\beta^0,j} W\right)^*
(K_{C_n}^*-1)\right\|_1\left((2\pi)^{-1}\int \parallel
f_{\beta^0}(\cdot)f_{X,Z}(x,\cdot)\parallel_1 dx\int_0^\tau
\eta_{\gamma^0}^2(t)dt\right)\\+\left\| \left(f^{(1)}_{\beta^0,j}f_{\beta^0} W\right)^*
    (K_{C_n}^*-1)\right\|_1\left((2\pi)^{-1}\int \parallel
f_{X,Z}(x,\cdot)\parallel_1 dx \int_0^\tau\eta_{\gamma^0}^2(t)dt\right).
 \end{multline*}
In the same way
$(1/2)\left\vert\mathbb{E}\left(\partial S_{n,1}(\theta)/\partial
\gamma_j\left.\right\vert_{\theta=\theta^0}\right)
%-\left({\partial S_{\theta^0,g}(\theta)}/{\partial
%\gamma_j}\left.\right\vert_{\theta=\theta^0}
%\right)
\right\vert$ is bounded by
\begin{multline*}
 \parallel
(f_{\beta^0} W)\star
    K_{C_n}-f_{\beta^0} W\parallel_\infty \left(\int \parallel f_{\beta^0}(\cdot)f_{X,Z}(x,\cdot)\parallel_1 dx\int_0^\tau
\eta^{(1)}_{\gamma^0,j}(t)\eta_{\gamma^0}(t)dt\right)\\+\parallel (f_{\beta^0}^2 W)\star
    K_{C_n}-f_{\beta^0}^2 W
\parallel_\infty \left(\int \parallel
f_{X,Z}(x,\cdot)\parallel_1 dx \int_0^\tau\eta^{(1)}_{\gamma^0,j}(t)\eta_{\gamma^0}(t)dt\right)\\
\leq  \parallel (f_{\beta^0} W)^*
(K_{C_n}^*-1)\parallel_1\left((2\pi)^{-1}\int \parallel
f_{\beta^0}(\cdot)f_{X,Z}(x,\cdot)\parallel_1 dx\int_0^\tau \eta^{(1)}_{\gamma^0,j}(t)\eta_{\gamma^0}(t)dt\right)\\+ \parallel (f_{\beta^0}^2 W)^*
    (K_{C_n}^*-1)\parallel_1\left((2\pi)^{-1}\int \parallel
f_{X,Z}(x,\cdot)\parallel_1 dx\int_0^\tau\eta^{(1)}_{\gamma^0,j}(t)\eta_{\gamma^0}(t)dt\right).
 \end{multline*}
Consequently
\begin{multline}
\label{b2b1}
\left\vert \mathbb{E}\left(\frac{\partial S_{n,1}(\theta)}{\partial
\beta_j}\left.\right\vert_{\theta=\theta^0}\right)
\right\vert \\
\leq
C(\theta^0)\left[\left\|\left(
f^{(1)}_{\beta^0,j}W\right)^*(K_{C_n}^*-1)\right\|_1+
\left\|\left(f^{(1)}_{\beta^0,j}f_{\beta^0}W
\right)^* (K_{C_n}^*-1)\right\|_1\right],
\end{multline}
and 
\begin{eqnarray}
\label{b2b2}
\qquad\left\vert \mathbb{E}\left(\frac{\partial S_{n,1}(\theta)}{\partial
\gamma_j}\left.\right\vert_{\theta=\theta^0}\right)
\right\vert
\leq C(\theta^0)\left[\left\|(f_{\beta^0}W)^*(K_{C_n}^*-1)\right\|_1 +\left\|(f_{\beta^0}^2W)^* (K_{C_n}^*-1)\right\|_1\right].
\end{eqnarray}
By combining \eref{b1b1}, \eref{b1b2}, \eref{b2b1} and \eref{b2b2} we get that
\begin{multline}
\label{bb1}
\left\vert \mathbb{E}\left(\frac{\partial S_{n,1}(\theta)}{\partial
\beta_j}\left.\right\vert_{\theta=\theta^0}\right)
\right\vert\leq C(\theta^0)
\\\times
\min\left\lbrace \left\|\left(f^{(1)}_{\beta^0,j}W\right)^*(K_{C_n}^*-1)\right\|_2^2
+\left\|\left(f^{(1)}_{\beta^0,j}f_{\beta^0}W\right)^*(K_{C_n}^*-1)\right\|_2^2,\right.\\\left.
\left\|\left(f^{(1)}_{\beta^0,j}W\right)^*(K_{C_n}^*-1)\right\|_1^2+
\left\|\left(f^{(1)}_{\beta^0,j}f_{\beta^0}W
\right)^* (K_{C_n}^*-1)\right\|_1^2\right\rbrace
\end{multline}
and
\begin{multline}
\label{bb2}
\left\vert \mathbb{E}\left(\frac{\partial S_{n,1}(\theta)}{\partial
\gamma_j}\left.\right\vert_{\theta=\theta^0}\right)
\right\vert\leq
C(\theta^0)\\\times\min\left\lbrace
\left\|(f_{\beta^0}W)^* (K_{C_n}^*-1)\right\|_2^2 +\left\|(f_{\beta^0}^2W)^*
(K_{C_n}^*-1)\right\|_2^2, \right.\\\left.
\left\|(f_{\beta^0}W)^* (K_{C_n}^*-1)\right\|_1^2+\left\|(f_{\beta^0}^2W)^* (K_{C_n}^*-1)\right\|_1^2\right\rbrace.
\end{multline}
According to Lemma \ref{contint}
\begin{eqnarray*}
\left\vert \mathbb{E}\left(\frac{\partial S_{n,1}(\theta)}{\partial
\beta_j}\left.\right\vert_{\theta=\theta^0}\right)
\right\vert^2=O\left(C_n^{-2a+1-r+(1-r)_-}\exp(-2dC_n^r)\right),
\end{eqnarray*}
and
\begin{eqnarray*}
\left\vert \mathbb{E}\left(\frac{\partial S_{n,1}(\theta)}{\partial
\gamma_j}\left.\right\vert_{\theta=\theta^0}\right)
\right\vert^2
=O\left(C_n^{-2a+1-r+(1-r)_-}\exp(-2dC_n^r)\right).
\end{eqnarray*}
\textbf{Study of the variance}

For the variance term, it is easy to see that 
\begin{eqnarray*}
\mbox{Var}\left(\frac{\partial
S_{n,1}(\theta)}{\partial \beta_j}\left.\right\vert_{\theta=\theta^0}\right)=\frac{8+o(1)}{n}\left[V_{1,j}+V_{2,j}\right],
\end{eqnarray*}
with
\begin{eqnarray*}
V_{1,j}=\mathbb{E} \left[\left(f^{(1)}_{\beta^0,j}f_{\beta^0}
W\right)\star
    K_{n,C_n}(U_1)\int_0^\tau \eta_{\gamma^0}^2(t)Y_1(t)dt\right]^2
\end{eqnarray*}
\mbox{ and }
\begin{eqnarray*}
V_{2,j}=\mathbb{E} \left[\left(f^{(1)}_{\beta^0,j}
W\right)\star
    K_{n,C_n}(U_1)\int_0^\tau \eta_{\gamma^0}(t)dN_1(t)\right]^2
\end{eqnarray*}
In the same way
\begin{eqnarray*}
\mbox{Var}\left(\frac{\partial
S_{n,1}(\theta)}{\partial \gamma_j}\right)=\frac{8+o(1)}{n}\left[V_{3,j}+V_{4,j}\right],
\end{eqnarray*}
with
\begin{eqnarray*}
V_{3,j}=\mathbb{E} \left[ (f_{\beta^0}^2
W)\star
    K_{n,C_n}(U_1)\int_0^\tau \eta^{(1)}_{\gamma^0,j}(t)\eta_{\gamma^0}(t)Y_1(t)dt\right]^2
\end{eqnarray*}
\mbox{ and }
\begin{eqnarray*}
V_{4,j}=\mathbb{E} \left[ (f_{\beta^0}
W)\star
    K_{n,C_n}(U_1)\int_0^\tau\eta^{(1)}_{\gamma^0,j}(t)dN_1(t)\right]^2.
\end{eqnarray*}
According to Lemma \ref{lemdeconv},
\begin{eqnarray*}
V_{1,j}&\leq&
\left[\int_0^\tau \eta_{\gamma^0}^2(t)dt \right]^2 \int \left\vert
\left<f_{X,Z}(x,\cdot)\star f_\varepsilon,\left(\left(f^{(1)}_{\beta^0,j}f_{\beta^0} W\right)\star K_{n,C_n}\right)^2 \right>\right\vert dx
\end{eqnarray*}
and
\begin{eqnarray*}
V_{3,j}&\leq&
\left[\int_0^\tau \eta^{(1)}_{\gamma^0,j}(t)\eta_{\gamma^0}(t)
dt \right]^2 \int \left\vert
\left<f_{X,Z}(x,\cdot)\star f_\varepsilon,\left((f_{\beta^0}^2 W)\star K_{n,C_n}\right)^2 \right>\right\vert dx.
\end{eqnarray*}
By applying the inequalities \eref{holder1} and \eref{holder2}
we get that
\begin{eqnarray*}
V_{1,j}&\leq&
\left[\int_0^\tau \eta_{\gamma^0}^2(t)dt \right]^2\\&&\times\min\left\lbrace \parallel f_\varepsilon\parallel_\infty
\left\|\left(f^{(1)}_{\beta^0,j}f_{\beta^0} W\right)^*\frac{K_{C_n}^*}{\overline{f_\varepsilon^*}} \right\|_2^2,
\left\|\left(f^{(1)}_{\beta^0,j}f_{\beta^0} W\right)^*\frac{K_{C_n}^*}{\overline{f_\varepsilon^*}} \right\|_1^2\right\rbrace,
\end{eqnarray*}
and
\begin{eqnarray*}
V_{3,j}&\leq&
\left[\int_0^\tau \eta^{(1)}_{\gamma^0,j}(t)\eta_{\gamma^0}(t)dt 
\right]^2\times\min\left\lbrace \parallel f_\varepsilon\parallel_\infty
\left\|\left(f_{\beta^0}^2 W\right)^*\frac{K_{C_n}^*}{\overline{f_\varepsilon^*}} \right\|_2^2,
\left\|\left(f_{\beta^0}^2 W\right)^*\frac{K_{C_n}^*}{\overline{f_\varepsilon^*}} \right\|_1^2\right\rbrace.
\end{eqnarray*}
Now, according to Lemma \ref{lemdeconv} we have
\begin{eqnarray*}
V_{2,j}&\leq&
\int \left\vert
\left<\varphi_{2}^2(x,\cdot)f_{X,Z}(x,\cdot)\star f_\varepsilon,
\left(\left( f^{(1)}_{\beta^0,j} W\right)\star K_{n,C_n}\right)^2 \right>\right\vert dx
\end{eqnarray*}
and
\begin{eqnarray*}
V_{4,j}&\leq&
\int \left\vert
\left<\varphi_{4,j}^2(x,\cdot)f_{X,Z}(x,\cdot)\star f_\varepsilon,\left(\left(f_{\beta^0}W\right)
\star K_{n,C_n}\right)^2 \right>\right\vert dx,
\end{eqnarray*}
where 
$$\varphi_{2}(X,Z)=\int_0^\tau \eta_{\gamma^0}(t)dN(t)\mbox{ and }\varphi_{4,j}(X,Z)=\int_0^\tau \eta^{(1)}_{\gamma^0,j}(t)dN(t).$$
By applying the inequalities \eref{holder1} and \eref{holder2}
we get that
\begin{eqnarray*}
V_{2,j}&\leq&
\mathbb{E}\left[\int_0^\tau \eta_{\gamma^0}(t)dN(t) \right]^2\\&&\times\min\left\lbrace \parallel f_\varepsilon\parallel_\infty
\left\|\left(f^{(1)}_{\beta^0,j} W\right)^*\frac{K_{C_n}^*}{\overline{f_\varepsilon^*}} \right\|_2^2,
\left\|\left(f^{(1)}_{\beta^0,j} W\right)^*\frac{K_{C_n}^*}{\overline{f_\varepsilon^*}} \right\|_1^2\right\rbrace,
\end{eqnarray*}
and
\begin{eqnarray*}
V_{4,j}&\leq&\mathbb{E}
\left[\int_0^\tau\eta^{(1)}_{\gamma^0,j}(t)dN(t) 
\right]^2\times\min\left\lbrace \parallel f_\varepsilon\parallel_\infty
\left\|\left(f_{\beta^0} W\right)^*\frac{K_{C_n}^*}{\overline{f_\varepsilon^*}} \right\|_2^2,
\left\|\left(f_{\beta^0} W\right)^*\frac{K_{C_n}^*}{\overline{f_\varepsilon^*}} \right\|_1^2\right\rbrace.
\end{eqnarray*}
The result follows by combining the bounds on the $V_{k,j}$'s for $k=1,\ldots,4$
and by applying Lemma \ref{contint} to get that
$$\mbox{Var}\left(\frac{\partial
S_{n,1}(\theta)}{\partial \gamma_j}\right)=O\left(C_n^{2(\alpha-a)+1-\rho+(1-\rho)_-}\exp(-2dC_n^r+2\delta
C_n^\rho))/n\right).$$

The proof of \textbf{3)} follows by choosing $C_n$, that realizes
the trade-off between the squared bias and the variance.

\begin{center}
\textbf{Proof of ii)}
\end{center}
According to \eref{Sn1}, $S_{n,1}^{(2)}(\theta^0)$ equals
\begin{eqnarray*}
\frac{\partial^2 S_{n,1}(\theta^0)}{\partial \theta^2}=\begin{pmatrix}
(S_{n,1}^{(2)})_{1,1} & (S_{n,1}^{(2)})_{1,2}\\
(S_{n,1}^{(2)})^\top_{1,2} & (S_{n,1}^{(2)})_{2,2}
\end{pmatrix},
\end{eqnarray*}
with 
\begin{eqnarray*}
(S_{n,1}^{(2)})_{1,2}&=&-\frac{2}{n}\sum_{i=1}^n(
f^{(1)}_{\beta^0}W)\star
K_{n,C_n}(U_i)\int_0^\tau\!\!\!(\eta_{\gamma^0}^{(1)}(t))^\top dN_i(t)\\&&+\frac{1}{n}\sum_{i=1}^n\left(f^{(1)}_{\beta^0}f_{\beta^0}W\right)\star
 K_{n,C_n}(U_i)\int_0^\tau(\eta^{(1)}_{\gamma^0}(t))^\top\eta_{\gamma^0}(t)Y_i(t)dt,
\end{eqnarray*}
\begin{eqnarray*}
(S_{n,1}^{(2)})_{1,1}&=&
-\frac{2}{n}\sum_{i=1}^n(
f^{(2)}_{\beta^0}W)\star
K_{n,C_n}(U_i)\int_0^\tau\!\!\eta_{\gamma^0}(t)dN_i(t)\\&&+\frac{1}{n}\sum_{i=1}^n\left(\frac{\partial^2
(f_{\beta}^2W)}{\partial \beta^2}\left.\right\vert_{\theta=\theta^0}\right)\star
 K_{n,C_n}(U_i)\int_0^\tau\!\!\!\eta_{\gamma^0}^2(t)Y_i(t)dt 
\end{eqnarray*}
and
\begin{eqnarray*}
(S_{n,1}^{(2)})_{2,2}(\theta)&=& -\frac{2}{n}\sum_{i=1}^n(
f_{\beta^0}W)\star
K_{n,C_n}(U_i)\int_0^\tau\!\!\eta_{\gamma^0}^{(2)}(t)dN_i(t)\\&&+\frac{1}{n}\sum_{i=1}^n(f_{\beta^0}^2W)\star
 K_{n,C_n}(U_i)\int_0^\tau\!\!\!\left(\frac{\partial^2 \eta_{\gamma}^2(t)}{\partial
\gamma^2}
\left.\right\vert_{\theta=\theta^0}\right)Y_i(t)dt.
\end{eqnarray*}
Under \textbf{\eref{super}}, for $C_n$ satisfying  \eref{condcons2},
$\mathbb{E}[S_{n,1}^{(2)}(\theta^0)-S_{\theta^0,g}^{(2)}(\theta^0)]^2=o(1)$. Hence \textbf{ii)} is proved.

\begin{center}
\textbf{Proof of iii)}
\end{center}
The proof of \textbf{iii)} follows by using the smoothness of $\beta\mapsto W f_\beta$ and
$\beta\mapsto Wf_\beta^2$ up to order 3, the smoothness of $\gamma\mapsto
\eta_\gamma$ and $\gamma\mapsto\eta_\gamma^2$ and by using the consistency of $\widehat{\theta}_1.$

\begin{center}
\textbf{Proof of iv)}
\end{center}
Let us introduce the random event  $E_n=\cap_{j,k}E_{n,j,k},$ where
$$
E_{n,j,k}=\left\lbrace \omega \mbox{ such that }\left\vert
\frac{\partial^2 S_{\theta^0,g}(\theta)}{\partial \theta_j\partial \theta_k}
\left.\right\vert_{\theta=\theta^0}
- \frac{\partial^2S_{n,1}(\theta,\omega)}{\partial \theta_j\partial \theta_k}
\left.\right\vert_{\theta=\theta^0} +(R_n)_{j,k}(\omega)\right\vert \leq
\frac{1}{2}\,\frac{\partial^2S_{\theta^0,g}(\theta)}{\partial \theta_j\partial \theta_k}
\left.\right\vert_{\theta=\theta^0} \right\rbrace.
$$
We first write that
\begin{eqnarray*}
\mathbb{E}\|\widehat{\theta}_1-\theta^0\|_{\ell^2}^2
&=&\mathbb{E}[\|\widehat{\theta}_1-\theta^0\|_{\ell^2}^2 \ind_{E_{n}}]
+\mathbb{E}[\|\widehat{\theta}_1-\theta^0\|_{\ell^2}^2 \ind_{E_n^c}]\\
&\leq&\mathbb{E}[\|\widehat{\theta}_1-\theta^0\|_{\ell^2}^2\ind_{E_n}]
+2\sup_{\theta\in \Theta}\|\theta\|_{\ell^2}^2\mathbb{P}(E_n^c).
\end{eqnarray*}
According to \eref{Rn} and \eref{base}  we have
\begin{eqnarray*}
\mathbb{E}[\|\widehat{\theta}_1-\theta^0\|_{\ell^2}^2\ind_{E_n}]
&\leq&\mathbb{E}\left[(S_{n,1}^{(1)}(\theta^0))^\top
[(S_{n,1}^{(2)}(\theta^0)+ R_n)^{-1}]^\top (S_{n,1}^{(2)}(\theta^0)+ R_n)^{-1}
S_{n,1}^{(1)}(\theta^0)\ind_{E_n}\right]\\
&\leq& C(m,p)\sup_{j,k}\left\vert \frac{\partial^2
S_{\theta^0,g}(\theta)}{\partial \theta_j\partial \theta_k}
\left.\right\vert_{\theta=\theta^0} \right\vert^{-2}
\mathbb{E}\left[(S_{n,1}^{(1)}(\theta^0))^\top S_{n,1}^{(1)}(\theta^0)\right]\\
&\leq& C(m,p)\sup_{j,k}\left\vert \frac{\partial^2S_{\theta^0,g}(\theta)}{\partial
      \theta_j\partial
      \theta_k}\left.\right\vert_{\theta=\theta^0} \right\vert^{-2}\varphi_n^2
\end{eqnarray*}
It remains thus to show that $\mathbb{P}(E_n^c)=o(\varphi^2_{n})$ with
$$\sup_{j,k}\mathbb{E}\left[\left(\frac{\partial^2(S_{n,1}(\theta)-S_{\theta^0,g}(\theta))}{\partial
      \theta_j\partial
      \theta_k}\left.\right\vert_{\theta=\theta^0}\right)^2\right]\leq \varphi^2_{n}.$$
We write
\begin{eqnarray*}
\mathbb{P}(E_n^c)&\leq&\sum_{j=1}^{m+p}\sum_{k=1}^{m+p}\mathbb{P}(E_{n,j,k}^c).
\end{eqnarray*} 
By Markov's inequality, for $q>2$,
\begin{eqnarray*}
\mathbb{P}(E_{n,j,k}^c)
&\leq&\left(\Big\vert \frac{1}{2}\,\frac{\partial^2S_{\theta^0,g}(\theta)}{\partial \theta_j\partial \theta_k}
\left.\right\vert_{\theta=\theta^0}\Big\vert^q\right)^{-1}\mathbb{E}\left[\left\vert \left(
\frac{\partial^2(S_{\theta^0,g}(\theta)-S_{n,1}(\theta))}{\partial \theta_j\partial \theta_k}
\left.\right\vert_{\theta=\theta^0}\right)
+(R_n)_{j,k} \right\vert^q   \right].
\end{eqnarray*}
In other words, using that $|a+b|^q\leq 2^{q-1}(|a|^q+|b|^q)$, we get
$$\left(\Big\vert \frac{1}{2}\,\frac{\partial^2S_{\theta^0,g}(\theta)}{\partial \theta_j\partial \theta_k}
\left.\right\vert_{\theta=\theta^0}\Big\vert^q\right)\mathbb{P}(E_{n,j,k}^c)$$
is less than
\begin{eqnarray*}
&& 2^{q-1}\left\vert  \left(
\frac{\partial^2S_{\theta^0,g}(\theta)}{\partial \theta_j\partial \theta_k}\left.\right\vert_{\theta=\theta^0}\right)
-\mathbb{E}\left[\left( \frac{\partial^2S_{n,1}(\theta)}{\partial \theta_j\partial \theta_k}
\left.\right\vert_{\theta=\theta^0}\right)\right]\right\vert^q\\&& +2^{q-1}\mathbb{E}\left[\left\vert
\mathbb{E}\left( \frac{\partial^2S_{n,1}(\theta)}{\partial \theta_j\partial \theta_k}
\left.\right\vert_{\theta=\theta^0}\right)
- \left(\frac{\partial^2S_{n,1}(\theta)}{\partial \theta_j\partial \theta_k}
\left.\right\vert_{\theta=\theta^0}\right)
+(R_n)_{j,k} \right\vert^q \right]\\
&\leq& 2^{q-1}\left\vert \left( \frac{\partial^2S_{\theta^0,g}(\theta)}{\partial \theta_j\partial \theta_k}
\left.\right\vert_{\theta=\theta^0}\right)
-\mathbb{E}\left[ \left(\frac{\partial^2 S_{n,1}(\theta)}{\partial \theta_j\partial \theta_k}
\left.\right\vert_{\theta=\theta^0}\right)\right]\right\vert^q\\
&&+2^{2q-2}
\left\lbrace\mathbb{E}\left[\left\vert
\mathbb{E}\left( \frac{\partial^2S_{n,1}(\theta)}{\partial \theta_j\partial \theta_k}
\left.\right\vert_{\theta=\theta^0}\right)
-\left( \frac{\partial^2S_{n,1}(\theta)}{\partial \theta_j\partial \theta_k}
\left.\right\vert_{\theta=\theta^0}\right)\right\vert ^q\right]
+\mathbb{E}\vert( R_n)_{j,k} \vert^q\right\rbrace.
\end{eqnarray*}
Now we apply the Rosenthal's inequality (see \eref{rosenthal} recalled in
Appendix),  to the sum of centered
variables
\begin{eqnarray*}\left(\frac{\partial^2S_{n,1}(\theta)}{\partial \theta_j\partial
    \theta_k}\left.\right\vert_{\theta=\theta^0}\right)
-\mathbb{E}\left[\left(\frac{\partial^2S_{n,1}(\theta)}{\partial \theta_j\partial
      \theta_k}\left.\right\vert_{\theta=\theta^0}\right)\right]&:=&n^{-1}\sum_{i=1}^nW_{n,i,j,k}.
\end{eqnarray*}
It follows that
\begin{multline*}
\mathbb{E}\left[\left\vert\left(\frac{\partial^2S_{n,1}(\theta)}{\partial \theta_j\partial
    \theta_k}\left.\right\vert_{\theta=\theta^0}\right)-\mathbb{E}\left(\frac{\partial^2S_{n,1}(\theta)}{\partial \theta_j\partial
    \theta_k}\left.\right\vert_{\theta=\theta^0}\right)\right\vert ^q\right]
\\\leq C(r)\left[n^{1-r}\mathbb{E}|W_{n,1,j,k}|^q+
n^{-q/2}\mathbb{E}^{q/2}|W_{n,1,j,k}|^2\right].
\end{multline*}
Take $q=4$ to get that
\begin{eqnarray*}
\mathbb{E}\left[\left\vert\left(\frac{\partial^2S_{n,1}(\theta)}{\partial \theta_j\partial
    \theta_k}\left.\right\vert_{\theta=\theta^0}\right)-\mathbb{E}\left(\left(\frac{\partial^2S_{n,1}(\theta)}{\partial \theta_j\partial
    \theta_k}\left.\right\vert_{\theta=\theta^0}\right)\right)\right\vert^4\right]\leq C(4)\left[n^{-3}\mathbb{E}|W_{n,1,j,k}|^4+
n^{-2}\mathbb{E}^{2}|W_{n,1,j,k}|^2\right].
\end{eqnarray*}
Therefore under the conditions ensuring that
% terme de biais+terme de variance traites en un seul coup.
\begin{eqnarray*}
\mathbb{E}\Big[\frac{\partial^2 S_{\theta^0,g}(\theta)}{\partial
\theta_j\theta_k}|_{\theta=\theta^0}-\frac{\partial^2
S_{n,1}(\theta)}{\theta_j\theta_k}|_{\theta=\theta^0}\Big]^2=o(1),\end{eqnarray*} we have
$$\mathbb{E}\Big[\frac{\partial^2 S_{\theta^0,g}(\theta)}{\partial
\theta_j\theta_k}|_{\theta=\theta^0}-\frac{\partial^2
S_{n,1}(\theta)}{\theta_j\theta_k}|_{\theta=\theta^0}\Big]^4=O(\varphi_{n}^4)=o(\varphi_n^2).$$
Now, by using the definition of $R_n$
and the smoothness properties of the derivatives
of $(W f_\beta)$ and $(W f_\beta^2)$ up to order 3, we get that
$\mathbb{E}((R_n)_{j,k}^4)=o(\|\widehat{\theta}_1-\theta^0\|_{\ell^2}^4)$, and we conclude that
\begin{eqnarray*}
\mathbb{E}\|\widehat{\theta}_1-\theta^0\|_{\ell^2}^2
&\leq& 4\mathbb{E}\left[(S_{n,1}^{(1)}(\theta^0))^\top
\Big[\Big(\frac{\partial^2 S_{\theta^0,g}(\theta)}{\partial
\theta_j\theta_k}|_{\theta=\theta^0}\Big)^{-1}\Big]^\top 
\Big(\frac{\partial^2 S_{\theta^0,g}(\theta)}{\partial
\theta_j\theta_k}|_{\theta=\theta^0}\Big)^{-1}
S_{n,1}^{(1)}(\theta^0)\right]
\\&&+o(\varphi_{n}^2)+o(\mathbb{E}[\|\widehat{\theta}_1-\theta^0\|_{\ell^2}^4]).\qquad
\qquad \Box
\end{eqnarray*}

\subsection{Proof of Theorem \ref{thCS}~: asymptotic normality}

According to Theorem \ref{thv1} and its proof, under \eref{C1}-\eref{C3}, $V_{n,j}(\theta^0)=O(1)$ and
the asymptotic normality of
$\widehat{\theta}_1$ follows by checking that
\\\textbf{v)} $\sqrt{n} \; S_{n,1}^{(1)}(\theta^0)\cvl \mathcal{N}\left(0, {\Sigma}_1\right)$, 
with ${\Sigma}_1$  defined in Theorem \ref{thCS}.

Let $H_{n,i}$, $\widehat{H}_{n,i}$,  $G_{n,i}$, and
$\widehat{G}_{n,i}$ be the processes defined 
for all $t \in [0,\tau]$ by 
\begin{eqnarray}
\label{hnic}
\qquad \widehat{H}_{n,i}(s)=\begin{pmatrix}
\frac{-2}{\sqrt{n}} 
(f_{\beta^0}^{(1)}W)\star
K_{n,C_n}(U_i)\eta_{\gamma^0}(s)
\\
\frac{-2}{\sqrt{n}} -(f_{\beta^0}W)\star
K_{n,C_n}(U_i)\eta_{\gamma^0}^{(1)}(s)
\end{pmatrix},
\end{eqnarray} 
\begin{eqnarray}
\label{hni}
\qquad {H}_{n,i}(s)=\begin{pmatrix}
\frac{-2}{\sqrt{n}} 
(f_{\beta^0}^{(1)}W)(Z_i)\eta_{\gamma^0}(s)
\\
\frac{-2}{\sqrt{n}} (f_{\beta^0}W)(Z_i)\eta_{\gamma^0}^{(1)}(s)
\end{pmatrix},
\end{eqnarray} 
\begin{eqnarray}
\label{gnic}
\widehat{G}_{n,i}(s)=
\begin{pmatrix}
\frac{2}{\sqrt{n}} 
(f_{\beta^0}f_{\beta^0}^{(1)}W)\star
K_{n,C_n}(U_i)\eta^2_{\gamma^0}(s)
\\
\frac{2}{\sqrt{n}} (f_{\beta^0}^2 W)\star
K_{n,C_n}(U_i)\eta_{\gamma^0}^{(1)}(s)\eta_{\gamma^0}(s)
\end{pmatrix},
\end{eqnarray}
\begin{eqnarray}
\label{gni}
\mbox{ and }\qquad{G}_{n,i}(s)=
\begin{pmatrix}
\frac{2}{\sqrt{n}} 
(f_{\beta^0}f_{\beta^0}^{(1)}W)(Z_i)\eta^2_{\gamma^0}(s)
\\
\frac{2}{\sqrt{n}} (f_{\beta^0}^2 W)(Z_i)\eta_{\gamma^0}^{(1)}(s)\eta_{\gamma^0}(s)
\end{pmatrix}.
\end{eqnarray}
According to \eref{Sn1b}, since $N_i(s)=M_i(s)+\Lambda_i(s,\theta^0,Z_i)$ (see
\eref{deflambda}), we get that
\begin{eqnarray*}
\sqrt{n} \; S_{n,1}^{(1)}(\theta^0)
&=&\sum_{i=1}^n\int_0^\tau
\widehat{H}_{n,i}(s)dN_i(s)+\sum_{i=1}^n\int_0^\tau
\widehat{G}_{n,i}(s)Y_i(s)ds\\
&=&A_1+A_2+A_3+A_4\end{eqnarray*}
with
\begin{eqnarray*}
A_1&=&\sum_{i=1}^n\int_0^\tau
H_{n,i}(s)dM_i(s),\quad A_2=\sum_{i=1}^n\int_0^\tau
[\widehat{H}_{n,i}(s)-H_{n,i}(s)]dM_i(s),\\
A_3&=&\sum_{i=1}^n\int_0^\tau
[\widehat{H}_{n,i}(s)-H_{n,i}(s)]d\Lambda_i(s,\theta^0,Z_i)
\mbox{ and }A_4=\sum_{i=1}^n\int_0^\tau
[\widehat{G}_{n,i}(s)-G_{n,i}(s)]Y_i(s)ds.
\end{eqnarray*}

\begin{center}
Study of $A_1$
\end{center}

The term $A_1$ is a linear combinations of stochastic integrals of
locally bounded and predictable processes, $H_{n,i}$, with respect to 
finite variation and local square integrable martingales, $M_i(\cdot)$. 
Consequently, $\mathbb{E}(A_1)=0$.
Denoting by  $< M >$ the predictable variation process of $ M $ 
 we have to verify the two following conditions 
for all $t$ in $[0,\tau]$ (see \cite{ABGK} page 68) :

\textbf{L1)} $\sum_{i=1}^n \int_0^t H_{n,i}(s)(H_{n,i}(s))^\top d<M_i>(s)\cvp
\widetilde{\Sigma}_1^2(t)$, 
with $\widetilde{\Sigma}_1^2(t)$ 
a positive covariance matrix defined by
\begin{eqnarray}
\tilde{\Sigma}^2_1(t)=4\mathbb{E}\left[\int_0^t
\begin{pmatrix}
(f_{\beta^0}^{(1)}W)(Z_i)\eta_{\gamma^0}(s)
\\
(f_{\beta^0}W)(Z_i)\eta_{\gamma^0}^{(1)}(s)
\end{pmatrix}
\begin{pmatrix}
(f_{\beta^0}^{(1)}W)(Z_i)\eta_{\gamma^0}(s)
\\
(f_{\beta^0}W)(Z_i)\eta_{\gamma^0}^{(1)}(s)
\end{pmatrix}^\top \eta_{\gamma^0}(s)Y_i(s)ds\right]
\end{eqnarray}

\textbf{L2)} For all $\epsilon >0$, $\sum_{i=1}^n \int_0^t H_{n,i}(s)(H_{n,i}(s))^\top\ind_{\| H_{n,i}(s)\|_{\ell^2}
\geq \epsilon}\;\;d<M_i>(s)=o_p(1)$ .
\\
Proof of \textbf{L1)}

Since  $<M_i>=\Lambda_i$, we have to prove that for all $t\in [0,\tau]$,
\begin{eqnarray}
\label{eq1}
\sum_{i=1}^n \int_0^t H_{n,i}(s)(H_{n,i}(s))^\top Y_i(s)f_{\beta^0}(Z_i)\eta_{\gamma^0}(s)ds\cvp 
\widetilde{\Sigma}_1^2(t).
\end{eqnarray}
We apply the following Lemma, which is a straightforward consequence of 
the fact that the set of functions $\mathcal{I}_t=\{x\mapsto \ind_{x\geq t}\}$ is a 
$\mathbb{P}$-Glivenko Cantelli class (see van der Vaart and Wellner \citeyear{VVWELLNER}).
\begin{lem}
For $j=1,\cdots,m$
\label{lGC1}
\begin{eqnarray*}
\sup_{0 \leq t \leq \tau}\left\vert  
\frac{1}{n}\sum_{i=1}^n Y_i(t)f_{\beta^0,j}(Z_i)(f_{\beta^0,j}^{(1)}W)(Z_i) - \mathbb{E}[Y(t)f_{\beta^0}(Z)(f_{\beta^0,j}^{(1)}W)(Z)]\right\vert\cvps 0,
\end{eqnarray*}
\begin{eqnarray*}
\sup_{0 \leq t \leq \tau}\left\vert  
\frac{1}{n}\sum_{i=1}^n Y_i(t)f_{\beta^0}^2(Z_i)W(Z_i) - \mathbb{E}[Y(t)f_{\beta^0}^2(Z)W(Z)]\right\vert\cvps 0
\end{eqnarray*}
\begin{eqnarray*}
\sup_{0 \leq t \leq \tau}\left\vert  
\frac{1}{n}\sum_{i=1}^n Y_i(t)f_{\beta^0}(Z_i)\vert(f_{\beta^0,j}^{(1)}W)(Z_i)\vert^3 - \mathbb{E}[Y(t)f_{\beta^0}(Z)\vert(f_{\beta^0,j}^{(1)}W)(Z)\vert^3]\right\vert\cvps 0,
\end{eqnarray*}
\begin{eqnarray*}
\mbox{ and }\sup_{0 \leq t \leq \tau}\left\vert  
\frac{1}{n}\sum_{i=1}^n Y_i(t)\vert f_{\beta^0}(Z_i)W(Z_i)\vert^3 -
\mathbb{E}[Y(t)\vert f_{\beta^0}(Z)W(Z)\vert^3]\right\vert\cvps 0.
\end{eqnarray*}
\end{lem}
Thus \textbf{L1)}
is checked .
\\
Proof of \textbf{L2)}. We have to check that for all $j=1,\ldots,m$
\begin{eqnarray*}
\frac{1}{n}\mathbb{E}\left[\sum_{i=1}^n\int_0^t \left[f_{\beta^0,j}^{(1)}(Z_i)W(Z_i)\eta_{\gamma^0}(s)\right]^2\ind_{\vert
f_{\beta^0,j}^{(1)}(Z_i)W(Z_i)\eta_{\gamma^0}(s)\vert \geq \epsilon\sqrt{n}}f_{\beta^0}(Z_i)\eta_{\gamma^0}(s)Y_i(s)ds\right]=o(1)
\end{eqnarray*}
and that for all $j=1,\ldots,p$
\begin{eqnarray*}
\frac{1}{n}\mathbb{E}\left[\sum_{i=1}^n\int_0^t \left[f_{\beta^0}^2(Z_i)W(Z_i)\eta^{(1)}_{\gamma^0,j}(s)\right]^2\ind_{\vert
f_{\beta^0}^2(Z_i)W(Z_i)\eta^{(1)}_{\gamma^0,j}(s)\vert \geq \epsilon\sqrt{n}}f_{\beta^0}(Z_i)\eta_{\gamma^0}(s)Y_i(s)ds\right]=o(1).
\end{eqnarray*}
This is a straightforward consequence of Lemma \ref{lGC} by writing that for
$j=1,\cdots,m $
\begin{multline*}
\frac{1}{n}\mathbb{E}\left[\sum_{i=1}^n\int_0^t \left[f_{\beta^0,j}^{(1)}(Z_i)W(Z_i)\eta_{\gamma^0}(s)\right]^2\ind_{\vert
f_{\beta^0,j}^{(1)}(Z_i)W(Z_i)\eta_{\gamma^0}(s)\vert \geq
\epsilon \sqrt{n}}f_{\beta^0}(Z_i)\eta_{\gamma^0}(s)Y_i(s)ds\right]\\
\leq
\frac{1}{n\sqrt{n}\epsilon}\mathbb{E}\left[\sum_{i=1}^n\int_0^t \vert f_{\beta^0,j}^{(1)}(Z_i)W(Z_i)\eta_{\gamma^0}(s)\vert^3f_{\beta^0}(Z_i)\eta_{\gamma^0}(s)Y_i(s)ds\right]
=o(1)
\end{multline*}
and for $j=1,\cdots,p$
\begin{multline*}
\frac{1}{n}\mathbb{E}\left[\sum_{i=1}^n\int_0^t \left[f_{\beta^0}^2(Z_i)W(Z_i)(\eta^{(1)}_{\gamma^0,j}(s)\right]^2\ind_{\vert
f_{\beta^0}^2(Z_i)W(Z_i)\eta^{(1)}_{\gamma^0,j}(s)\vert \geq
\epsilon\sqrt{n}}f_{\beta^0}(Z_i)\eta_{\gamma^0}(s)Y_i(s)ds\right]\\
\leq\frac{1}{\epsilon n\sqrt{n}}\mathbb{E}\left[\sum_{i=1}^n\int_0^t
\vert f_{\beta^0}^2(Z_i)W(Z_i)\vert^3\vert \eta^{(1)}_{\gamma^0,j}(s)\vert ^3 f_{\beta^0}(Z_i)\eta_{\gamma^0}(s)Y_i(s)ds\right]
=o(1).
\end{multline*}
Thus \textbf{L2)} is checked.

\begin{center}
Study of $A_2$
\end{center}
Since 
$\mathbb{E}(A_2)=0$, we use the following lemma, analogous to Lemma \ref{lGC1}.
\begin{lem}
\label{lGC}
Under \textbf{\eref{condfeps}}-\textbf{\eref{super}}, for $C_n$ satisfying
\eref{condcons2} thenfor $j=1,\cdots, m$

\begin{eqnarray*}
\sup_{0 \leq t \leq \tau}\left\vert  
\frac{1}{n}\sum_{i=1}^n Y_i(t)f_{\beta^0}(Z_i)(f_{\beta^0,j}^{(1)}W)\star
K_{n,C_n}(U_i) - \mathbb{E}[Y(t)f_{\beta^0}(Z)(f_{\beta^0,j}^{(1)}W)(Z)]\right\vert\cvps 0,
\end{eqnarray*}
\begin{eqnarray*}
\mbox{ and }\sup_{0 \leq t \leq \tau}\left\vert  
\frac{1}{n}\sum_{i=1}^n Y_i(t)f_{\beta^0}(Z_i)(f_{\beta^0}W)\star
K_{n,C_n}(U_i) - \mathbb{E}[Y(t)f_{\beta^0}^2(Z)W(Z)]\right\vert\cvps 0.
\end{eqnarray*}
\end{lem}
It follows that $A_2=o_p(1).$

\begin{center}
Study of $A_3$
\end{center}

It is noteworthy that the term $A_3$  can be seen as
triangular arrays of row-wise independent centered random variables that is
$$A_3=\sum_{i=1}^n V_{n,i}+\mathbb{E}(A_3),$$
with $\sum_{i=1}^n V_{n,i}=A_3-\mathbb{E}(A_3)$.
Consequently, the asymptotic normality  follows
by checking that

\textbf{v-a)} $\mathbb{E}(A_3)=o_p(1)$

\textbf{v-b)} $\sum_{i=1}^n\mathbb{E}[(V_{n,i})^2]\cv \Sigma_3^2$

\textbf{v-c)} For all $\epsilon>0$, $\sum_{i=1}^n\mathbb{E}[(V_{n,i})^2\ind_{\parallel
V_{n,i}\parallel_{\ell^2}\geq \epsilon}]\cv 0$ (Lindeberg Condition).

By definition, $A_3$ equals
\begin{eqnarray*}
-\frac{2}{\sqrt{n}}\sum_{i=1}^n\int_0^\tau\begin{pmatrix}(f_{\beta^0}^{(1)}W)\star
K_{n,C_n}(U_i)-(f_{\beta^0}^{(1)}W)(Z_i)\eta_{\gamma^0}(s)\\
 (f_{\beta^0}W)\star
K_{n,C_n}(U_i)-(f_{\beta^0}W)(Z_i) \eta^{(1)}_{\gamma^0}(s) \end{pmatrix}
Y_i(s)f_{\beta^0}(Z_i)\eta_{\gamma^0}(s)ds.
\end{eqnarray*}

Let us start with the study of the variance (\textbf{v-b}). Under \eref{C1}-\eref{C3}
\begin{eqnarray*}
\mbox{Var}\left[-\frac{2}{\sqrt{n}}\sum_{i=1}^n\int_0^\tau(f_{\beta^0}^{(1)}W)\star
K_{n,C_n}(U_i)-(f_{\beta^0}^{(1)}W)(Z_i)\eta_{\gamma^0}(s)
Y_i(s)f_{\beta^0}(Z_i)\eta_{\gamma^0}(s)ds\right]=O(1),
\end{eqnarray*}
and
\begin{eqnarray*}
\mbox{Var}\left[-\frac{2}{\sqrt{n}}\sum_{i=1}^n\int_0^\tau
 (f_{\beta^0}^{(1)}W)\star
K_{n,C_n}(U_i)-(f_{\beta^0}W)(Z_i) \eta^{(1)}_{\gamma^0}(s)
Y_i(s)f_{\beta^0}(Z_i)\eta_{\gamma^0}(s)ds
\right]=O(1).
\end{eqnarray*}
It follows that \textbf{v-b)} is checked.

We now come to the bias term and write that
\begin{eqnarray*}
\mathbb{E}(A_3)\!\!\!&=&\!\!\!
\mathbb{E}\left\lbrace-\frac{2}{\sqrt{n}}\sum_{i=1}^n\int_0^\tau\begin{pmatrix}(f_{\beta^0}^{(1)}W)\star
K_{n,C_n}(U_i)-(f_{\beta^0}^{(1)}W)(Z_i)\eta_{\gamma^0}(s)\\
 (f_{\beta^0}W)\star
K_{n,C_n}(U_i)-(f_{\beta^0}W)(Z_i) \eta^{(1)}_{\gamma^0}(s) \end{pmatrix}
Y_i(s)f_{\beta^0}(Z_i)\eta_{\gamma^0}(s)   \right\rbrace ds\\
\!\!\!&=&\!\!\!-2\sqrt{n}\begin{pmatrix} 
\mathbb{E}\left\lbrace \left[(f_{\beta^0}^{(1)}W)\star
K_{n,C_n}(U)-(f_{\beta^0}^{(1)}W)(Z)\right] f_{\beta^0}(Z_i)\int_0^\tau\eta_{\gamma^0}^2(s)Y(s)
ds
\right\rbrace \\
\mathbb{E}\left\lbrace \left[(f_{\beta^0}W)\star
K_{n,C_n}(U_i)-(f_{\beta^0}W)(Z_i)  \right] f_{\beta^0}(Z_i)\int_0^\tau\eta^{(1)}_{\gamma^0}(s)\eta_{\gamma^0}(s)
Y(s) ds\right\rbrace \end{pmatrix}.
\end{eqnarray*}
According to Lemma \ref{lemdeconv}
\begin{eqnarray*}
\mathbb{E}(A_3)
\!\!\!&=&\!\!\!-2\sqrt{n}\begin{pmatrix} 
\mathbb{E}\left\lbrace\left[(f_{\beta^0}^{(1)}W)\star
K_{C_n}(Z)-(f_{\beta^0}^{(1)}W)(Z)\right]f_{\beta^0}(Z)
 \int_0^\tau Y_(s)\eta_{\gamma^0}^2(s)ds \right\rbrace\\~\\
\mathbb{E}\left\lbrace \left[(f_{\beta^0}W)\star
K_{C_n}(Z)-(f_{\beta^0}W)(Z) 
\right] f_{\beta^0}(Z) \int_0^\tau Y(s)\eta^{(1)}_{\gamma^0}(s)\eta_{\gamma^0}(s)ds\right\rbrace\end{pmatrix}\\~\\
\!\!\!&=&\!\!\!-2\sqrt{n}\begin{pmatrix} 
\mathbb{E}\left\lbrace\left[(f_{\beta^0}^{(1)}W)\star
K_{C_n}(Z)-(f_{\beta^0}^{(1)}W)(Z)\right]f_{\beta^0}(Z)
 \int_0^\tau Y(s)\eta_{\gamma^0}^2(s)ds \right\rbrace\\~\\
\mathbb{E}\left\lbrace \left[(f_{\beta^0}W)\star
K_{C_n}(Z)-(f_{\beta^0}W)(Z)  \right]
f_{\beta^0}(Z)
\int_0^\tau Y(s )\eta^{(1)}_{\gamma^0}(s)\eta_{\gamma^0}(s)ds\right\rbrace\end{pmatrix}\\~\\
\!\!\!&=&\!\!\!-2\displaystyle\sqrt{n}\begin{pmatrix} 
\int \left<(f_{\beta^0}^{(1)}W)\star
K_{C_n}(z)-(f_{\beta^0}^{(1)}W)(z),f_{\beta^0}(z)f_{X,Z}(x,z)\right>
\left( \int_0^\tau\ind_{x\geq s}\,
\eta_{\gamma^0}^2(s)ds\right) \,dx \\~\\
\int \left<(f_{\beta^0}W)\star
K_{C_n}(z)-(f_{\beta^0}W)(z),f_{\beta^0}(z)\right>
\left(\int_0^\tau\ind_{x\geq s}
\eta^{(1)}_{\gamma^0}(s)\eta_{\gamma^0}(s)ds\right)dx\end{pmatrix}.
\end{eqnarray*}
For $j=1,\ldots,m$
\begin{eqnarray*}
\left\vert \int \left<(f_{\beta^0,j}^{(1)}W)\star
K_{C_n}(z)-(f_{\beta^0,j}^{(1)}W)(z),f_{\beta^0}(z)f_{X,Z}(x,z)\right>
 \left(\int_0^\tau\ind_{x\geq s}
\eta_{\gamma^0}^2(s)ds\right)dx\right\vert
\end{eqnarray*}
is less than
\begin{eqnarray*}
\left(\int_0^\tau\eta_{\gamma^0}^2(s)ds\right)\int \left\vert \left<(f_{\beta^0,j}^{(1)}W)\star
K_{C_n}(z)-(f_{\beta^0,j}^{(1)}W)(z),f_{\beta^0}(z)f_{X,Z}(x,z)\right>\right\vert
dx,
\end{eqnarray*}
which is, according to \eref{holder1} and \eref{holder2},
less than
\begin{multline*}
\left(\int_0^\tau\eta_{\gamma^0}^2(s)ds\right) \min\left\lbrace\int \left\|(f_{\beta^0,j}^{(1)}W)\star
K_{C_n}-(f_{\beta^0,j}^{(1)}W)\right\|_2\left\|f_{\beta^0}(\cdot)f_{X,Z}(x,\cdot)\right\|_2
dx 
, \right.\\\left.\int \left\|(f_{\beta^0,j}^{(1)}W)\star
K_{C_n}-(f_{\beta^0,j}^{(1)}W)\right\|_\infty \left\|f_{\beta^0}(\cdot)f_{X,Z}(x,\cdot)\right\|_1
dx \right\rbrace 
\end{multline*}
that is less than 
\begin{multline*}
(2\pi)^{-1}\left(\int_0^\tau\eta_{\gamma^0}^2(s)ds\right)\times\min\left\lbrace \left\|(f_{\beta^0,j}^{(1)}W)^*(K_{C_n}^*-1)\right\|_2\int\left\|f_{\beta^0}(\cdot)f_{X,Z}(x,\cdot)\right\|_2
dx 
, \right.\\\left. \left\|(f_{\beta^0,j}^{(1)}W)^*(K_{C_n}^*-1)\right\|_\infty \int\left\|f_{\beta^0}(\cdot)f_{X,Z}(x,\cdot)\right\|_1
dx \right\rbrace .
\end{multline*}
In the same way we obtain that for $j=1,\ldots,p$
\begin{eqnarray*}
\left\vert\int \left<(f_{\beta^0}W)\star
K_{C_n}(z)-(f_{\beta^0}W)(z),f_{\beta^0}(z)\right>
\left(\int_0^\tau\ind_{x\geq s}
\eta^{(1)}_{\gamma^0}(s)\eta_{\gamma^0}(s)ds\right)dx\right\vert
\end{eqnarray*}
is less than
\begin{multline*}
(2\pi)^{-1}\left(\int_0^\tau\vert
\eta_{\gamma^0}(s)\eta_{\gamma^0,j}^{(1)}(s)\vert ds\right)\times\min\left\lbrace \left\|(f_{\beta^0}W)^*(K_{C_n}^*-1)\right\|_2\int\left\|f_{\beta^0}(\cdot)f_{X,Z}(x,\cdot)\right\|_2
dx 
, \right.\\\left. \left\|(f_{\beta^0}W)^*(K_{C_n}^*-1)\right\|_\infty \int\left\|f_{\beta^0}(\cdot)f_{X,Z}(x,\cdot)\right\|_1
dx \right\rbrace.
\end{multline*}
Consequently, under \textbf{\eref{super}},
$\mathbb{E}(A_3)=O(\sqrt{n}C_n^{-a+(1-r)/2+(1-r)_-/2\exp(-dC_n^r)}).$
Under \eref{C1}-\eref{C3}, $\mbox{Var}(A_3)=O(1)$
and hence $C_n$ can be chosen such that $\mathbb{E}(A_3)=o(1).$
It follows that \textbf{v-a)} is checked.

In order to check the Lindeberg condition we write that for $j=1,\cdots,m$
\begin{multline*}
\frac{1}{n}\mathbb{E}\left[\sum_{i=1}^n\int_0^t \left[(f_{\beta^0,j}^{(1)}W)\star K_{n,C_n}(U_i)\eta_{\gamma^0}(s)\right]^2\ind_{\vert
(f_{\beta^0,j}^{(1)}W)\star K_{n,C_n}(U_i) \eta_{\gamma^0}(s)\vert \geq
\epsilon \sqrt{n}}f_{\beta^0}(Z_i)\eta_{\gamma^0}(s)Y_i(s)ds\right]\\
\leq
\frac{1}{n\sqrt{n}\epsilon}\mathbb{E}\left[\sum_{i=1}^n\int_0^t \vert (f_{\beta^0,j}^{(1)}W)\star K_{n,C_n}(U_i)\eta_{\gamma^0}(s)\vert^3f_{\beta^0}(Z_i)\eta_{\gamma^0}(s)Y_i(s)ds\right]
=o(1)
\end{multline*}
and for $j=1,\cdots,p$
\begin{multline*}
\frac{1}{n}\mathbb{E}\left[\sum_{i=1}^n\int_0^t \left[(f_{\beta^0}^2W)\star K_{n,C_n}(U_i)\eta^{(1)}_{\gamma^0,j}(s)\right]^2\ind_{\vert
(f_{\beta^0}^2W)\star K_{n,C_n}(U_i)\eta^{(1)}_{\gamma^0,j}(s)\vert \geq
\epsilon\sqrt{n}}f_{\beta^0}(Z_i)\eta_{\gamma^0}(s)Y_i(s)ds\right]\\
\leq\frac{1}{\epsilon n\sqrt{n}}\mathbb{E}\left[\sum_{i=1}^n\int_0^t
\vert (f_{\beta^0}^2W)\star K_{n,C_n}(U_i)\vert^3\vert \eta^{(1)}_{\gamma^0,j}(s)\vert ^3 f_{\beta^0}(Z_i)\eta_{\gamma^0}(s)Y_i(s)ds\right]
=o(1).
\end{multline*}
It follows that \textbf{v-c)} is checked.
\begin{center}
Study of $A_4$
\end{center}
The study of $A_4$, quite similar to the study of $A_3$ is omitted. \hfill $\Box$

\subsection{Proof of Theorem \ref{thv2}~:}
The proof of Thorem \ref{thv2}, quite classical is omitted.
\section{Appendix}
\setcounter{equation}{0}
\setcounter{lem}{0}
\setcounter{theo}{0}

\begin{lem}
\label{lemdeconv}
Let $\varphi$ be such that $\mathbb{E}(\vert \varphi(X,Z)\vert)$ is finite and
let $\Phi$ such that $\mathbb{E}(\vert \Phi(U)\vert)$ is finite.
Under the assumptions \eref{hyp3} and \eref{hyp2}, then
\begin{eqnarray*}
\mathbb{E}[\varphi(X,Z)\Phi\star
K_{n,C_n}(U)]=\mathbb{E}[\varphi(X,Z)\Phi\star K_{C_n}(Z)],
\end{eqnarray*}
and
\begin{eqnarray*}
\mathbb{E}\left[\varphi(X,Z)\Phi\star K_{n,C_n}(U)\right]^2=\int
\left<\left(\varphi^2(x,\cdot)f_{X,Z}(x,\cdot)\right)\star f_\varepsilon,(\Phi\star K_{n,C_n})^2\right>dx.
\end{eqnarray*}
\end{lem}

\textbf{Proof of Lemma \ref{lemdeconv}~:}
Set $f_{X,U,Z}$
the joint distribution of $(X,U,Z)$. Under \eref{hyp3} and
\eref{hyp2}, $f_{X,U,Z}(x,u,z)=f_{X,Z}(x,z)f_\varepsilon(u-z)$.
Hence by the
Parseval's formula
\begin{eqnarray*}
\mathbb{E}\left[\varphi(X,Z)\Phi\star
  K_{n,C_n}(U)\right]
%&=&\iiint \varphi(x,z)\Phi\star K_{n,C_n}(u)f_{X,U,Z}(x,u,z)dx\,du\,dz\\
&=&\iiint \varphi(x,z)\Phi\star K_{n,C_n}(u)f_{X,Z}(x,z)f_\varepsilon(u-z)du\,dx\,dz\\
&=&\iint \varphi(x,z)f_{X,Z}(x,z)\int\Phi\star
K_{n,C_n}(u)f_\varepsilon(u-z)du\,dx\,dz\\
&=&(2\pi)^{-1}\iint \varphi(x,z)f_{X,Z}(x,z)\int\Phi^*(y)K_{n,C_n}^*(y)f_\varepsilon^*(y)e^{-iyz}dy\,dx\,dz\\
&=&(2\pi)^{-1}\iint
\varphi(x,z)f_{X,Z}(x,z)\int\Phi^*(y)\frac{K_{C_n}^*(y)}{
\overline{f_\varepsilon^*}(y)}\overline{f_\varepsilon^*}(y)e^{-iyz}dy\,dx\,dz\\
&=&(2\pi)^{-1}\iint
\varphi(x,z)f_{X,Z}(x,z)\int\Phi^*(y)K_{C_n}^*(y)e^{-iyz}dy\,dx\,dz\\
&=&\iint
\varphi(x,z)f_{X,Z}(x,z)\int\Phi(u)K_{C_n}(z-u)du\,dx\,dz\\
&=&\iint
\varphi(x,z)\Phi\star K_{C_n}(z)f_{X,Z}(x,z)\,dx\,dz.
\end{eqnarray*}
In the same way,
\begin{eqnarray*}
\mathbb{E}\left[\varphi(X,Z)\Phi\star K_{n,C_n}(U)\right]^2
%&=& \iiint \varphi^2(x,z)(\Phi\star K_{n,C_n}(u))^2
%f_{X,U,Z}(x,u,z)dx\,du\,dz\\
&=&
\iiint \varphi^2(x,z)(\Phi\star K_{n,C_n}(u))^2 f_{X,Z}(x,z)f_\varepsilon(u-z)dx\,du\,dz\\&=&
\iiint \varphi^2(x,z)(\Phi\star K_{n,C_n}(u))^2
f_{X,Z}(x,z)f_\varepsilon(u-z)dx\,du\,dz\\&=&\int \left< (\varphi^2(x,\cdot)f_{X,Z}(x,\cdot))\star
f_\varepsilon,(\Phi\star K_{n,C_n})^2\right>dx.\qquad\qquad \Box
\end{eqnarray*}

\begin{lem}
\label{contint}
For $a$, $r$ two nonnegative numbers,
Then
\begin{eqnarray}
\label{bexp}
\int_{\vert u\vert \geq C_n}\vert u\vert^{-\nu}\exp(-\lambda\vert
u\vert^\delta) du\leq \frac{1}{C(\nu,\lambda,\delta)}C_n^{-\nu+1-\delta}\exp\{-\lambda C_n^\delta\}.
\end{eqnarray}
Furthermore, if $f_{\varepsilon}$ satisfies \textbf{\eref{condfeps}}, then
\begin{eqnarray*}
\int_{\vert u\vert \leq C_n}\frac{\vert u\vert^{-\nu}\exp(-\lambda\vert u\vert^\delta)}{\vert
  f_{\varepsilon}^*(u)\vert}du \leq
\frac{1}{C(\alpha,\delta,\rho,\nu,\lambda,\delta)\underline{C}(f_{\varepsilon})}\max[1,C_n^{(\alpha-\nu+1-\delta)}\exp\{-\lambda
C_n^\delta+\delta
C_n^\rho\}].
\end{eqnarray*}
\end{lem}

\begin{lem}
\noindent \textbf{Rosenthal's inequality }(Rosenthal~(1970),
Petrov~(1995)).
For $U_1,\ldots,U_n$, be $n$ independent centered
random variables, there exists a constant $C(r)$ such that for
$r\geq 1$,
\begin{eqnarray}
\label{rosenthal}
{\mathbb E}[|\sum_{i=1}^n U_i|^r] \leq
C(r) [ \sum_{i=1}^n {\mathbb E}[|U_i|^r] +( \sum_{i=1}^n
{\mathbb E}[U_i^2])^{r/2}].
\end{eqnarray}\end{lem}

\bibliographystyle{plain}
\bibliography{Biblio}
%\bibliography{/home/mlm/DOC/PAPER/A2/biblio}
{\small
\parbox{8cm}{{\sc Marie-Laure MARTIN-MAGNIETTE}\\
{\sc Institut National Agronomique Paris-Grignon\\
 Math\'ematique et Informatique
Appliqu\'ees, \\
16, rue Claude Bernard\\
75231 Paris cedex 05, France,}\\~\\
{\sc Institut National de la Recherche Agronomique, \\
Unit\'e de Recherche en G\'enomique V\'eg\'etale,\\
UMR INRA 1165- CNRS 8114 -UEVE\\
2,rue Gaston Cr\'emieux- CP 5708\\
91057 Evry Cedex, France} \\
e-mail~: mlmartin@inapg.fr}
\hfill \parbox{8cm}{{\sc Marie-Luce TAUPIN}\\
{\sc Laboratoire de Probabilit\'es, Statistique et Mod\'elisation, UMR 8628,}\\
{\sc Universit\'e Paris-Sud, B\^at. 425,}\\
{\sc 91405 Orsay Cedex, France}\\
e-mail~: marie-luce.taupin@math.u-psud.fr}}

\end{document}